\documentclass[11pt]{amsart}
\usepackage{amssymb,amsmath,amsthm,andy}

\newtheorem{thm}{Theorem}[section]
\newtheorem{prop}[thm]{Proposition}
\newtheorem{lem}[thm]{Lemma}
\newtheorem{cor}[thm]{Corollary}

\newtheorem{rem}[thm]{Remark}
\newcommand{\Null}[1]{\text{null}(#1)}

\newcommand{\Vtp}{V_{\tau p}}

\newcommand{\tnabla}{\widetilde\nabla}

\begin{document} 

\begin{abstract}Let $p:\C\to\R$ be a subharmonic, nonharmonic polynomial and $\tau\in\R$ a parameter.
Define $\Zbstp = \frac{\p}{\p\z} + \tau\frac{\p p}{\p \z} = e^{-\tau p}\frac{\p}{\p\z} e^{\tau p}$, 
a closed, densely defined operator on 
$L^2(\C)$. If $\Boxtp = \Zbstp\Zbstp^*$ and $\Boxwtp = \Zbstp^*\Zbstp$, we solve
the heat equations $\p_s u + \Boxtp u=0$, $u(0,z)=f(z)$ and $\p_s \tilde u + \Boxwtp \tilde u=0$, $\tilde u(0,z) 
= \tilde f(z)$. We write the solutions via heat semigroups and show that the solutions can be written as
integrals against distributional kernels. We prove that the kernels are $C^\infty$ off of the diagonal
$\{(s,z,w) : s=0 \text{ and } z=w\}$ and find pointwise bounds for the kernels and their derivatives.
\end{abstract}

\title [Pointwise Estimates of Heat Kernels in $\R\times\C$] 
{Pointwise Estimates for Relative Fundamental Solutions of Heat Equations in $\R\times\C$}
\author{Andrew Raich}

\address{
Department of Mathematics\\ Texas A\&M University\\ Mailstop 3368  \\ College Station, TX
77843-3368}



\maketitle
%
%
\section{Introduction}\label{sec:intro}

The object of this article is to study the relative fundamental solutions of a class of
heat equations on $\R\times\C$ that are motivated by and have applications to questions in
several complex variables. We solve each heat equation via a heat semigroup and write the solution
as a fractional integral operator. We find the regularity of the integral kernel and find pointwise
estimates on the kernel and its derivatives. One additional point of interest is that the 
infintesimal generator of the semigroup is a  magnetic Schr\"odinger operator whose
electric potential is nonpositive, yet the large time behavior of the semigroup is well-controlled.

Let $p:\C\to\R$ be a subharmonic, nonharmonic polynomial and $\tau\in\R$ a parameter. If $z = x_1 + i x_2$ and
$\frac{\p}{\p\z} = \frac 12 \left(\frac{\p}{\p x_1} + i \frac{\p}{\p x_2}\right)$, define $\Zbstp$
to be the operator 
\[
\Zbstp  = \frac{\p}{\p\z} +  \tau \frac{\p p}{\p\z},
\] 
and let $\Zstp  = -\Zbstp^* = \frac{\p}{\p z} - \tau \frac{\p p}{\p z}$ be the negative of the
$L^2$-adjoint of $\Zbstp $. If  $\Boxtp = -\Zbstp  \Zstp $ and $\Boxwtp = -\Zstp\Zbstp$, 
then our goal is to  understand the
heat equations:
\begin{equation}\label{eq:he i1}
\begin{cases} {\displaystyle \frac{\p u}{\p s}} + \Boxtp u=0\vspace*{.1in}\\
u(0,z) = f(z)\end{cases}
\end{equation}
and
\begin{equation}\label{eq:he2}
\begin{cases} {\displaystyle \frac{\p \tilde u}{\p s}} + \Boxwtp \tilde u=0\vspace*{.1in}\\
\tilde u(0,z) = \tilde f(z).\end{cases}
\end{equation}	
We  write our solutions 
\begin{align}
u(s,z) &= e^{-s\Boxtp}[f](z) = \int_{\C}\Htp(s,z,w)f(w)\, dA(w) \label{eq:he I}\\
\intertext{and}
\tilde u(s,z) &= e^{-s\Boxwtp}[f](z) = \int_{\C}\Hwtp(s,z,w)f(w)\, dA(w) \label{eq:hew I}
\end{align}
where $dA$ is Lebesgue meausure on $\C$.
After determining the
regularity of the kernels $\Htp(s,z,w)$ and $\Hwtp(s,z,w)$, we obtain pointwise estimates 
on the two functions and their derivatives.

\subsection{Background}\label{subsec:background}
Since $\Zbstp = e^{-\tau p}\frac{\p}{\p\z}e^{\tau p}$, 
solving $\Zbstp u=f$ in $L^2(\C)$ is equivalent to solving $\dbar u=f$ in $L^2(\C,e^{-2\tau p})$. 
When $\tau=1$, the 
$\ZZ_p$-problem has been solved in $L^2(\C)$ \cite{Christ91} and on smoothly bounded domains in $\C$
\cite{Ber96}. Christ proves that $G_p = \Box_p^{-1}$ is a well-defined, bounded operator on $L^2(\C)$, so
$R_p = Z_p G_p$ is 
the relative fundamental solution of $\ZZ_p$. Christ writes $G_p$ and $R_p$ as integral operators
and finds pointwise bounds on the integral kernels. Christ's techniques involve establishing the
$L^2$-theory of $\ZZ_p$ and using integration by parts to establish the decay. Christ proves that the
kernels of $ G_p$ and $R_p$ decay exponentially with decay term $e^{-\rho(z,w)}$ where $\rho$
is a metric whose size is governed by the size of the derivatives of $p$. 
Berndtsson also
analyzes $G_p$ and $R_p$ except that he recognizes that $\Box_p$ is a 
magnetic Schr\"odinger operator
of the form  $2\Box_p = \frac 12(i\nabla-a)^2 +V$ where $a = (-\frac{\p p}{\p x_2},
\frac{\p p}{\p x_1})$ and $V = \frac 12\triangle p$. 

The weighted $\dbar$-problem in $\C$ is not the only problem in which $\Boxtp$ is the crucial
object to understand.
Mathematicians have analyzed operators on Hartogs domains in $\C^n$ by understanding
weighted operators on the base space. The original operators are then reconstructed by 
Fourier series \cite{Li89, FoSi91, Ber94}.
Recently, on a class of Hartogs domains $\Omega\subset\C^2$,
Fu and Straube \cite{FuSt02,FuSt04} establish an equivalence between the compactness of
the $\dbar$-Neumann problem and the blowup of the smallest eigenvalue of
$\Boxtp$ as $\tau\to\infty$. Christ and Fu \cite{ChFu05} build on 
the work of Fu and Straube to show that the following
are equivalent: compactness of the inverse operator of $\dbar$-Neumann Laplacian, compactness of
the inverse operator of Kohn Laplacian $\Box_b$, and $b\Omega$ satisfiying property $(P)$.

In \cite{Rai05h}, we present an alternative method to analyzing $\Boxtp$. 
We study (\ref{eq:he i1})
in the case $\tau>0$ (and $\tau=1$ in particular) and show that (\ref{eq:he I}) solves (\ref{eq:he i1}).
We prove that $\int_0^\infty e^{-s\Boxtp}\, ds = \Boxtp^{-1}$ is well-defined and find pointwise
estimates of $\Htp(s,z,w)$ and its derivatives. The techniques include the spectral theorem, results
about one-parameter families (OPF) 
of operators \cite{Rai05f}, the Feynman-Kac-It\^o formula from mathematical
physics, and an energy inequality. Our results allow us to recover estimates of $G_p$ and $R_p$. In the
cases for which $\rho$ has 
been computed, the decay in \cite{Christ91} and \cite{Rai05h} agree.

The motivation of \cite{Rai05h, Rai05f} is not, however, to reprove Christ's results but rather to use the
strategy of Nagel and Stein \cite{NaSt00, NaSt04, NaSt03} to understand $e^{-s\Box_b}$ 
on (the boundary of) decoupled domains. 
The boundary of a decoupled domain is a manifold
$M = \{(z_1, \dots, z_n) \in \C^n : \Imm z_n = P(z_1,\dots, z_{n-1})\}$ where $P(z_1,\dots,z_{n-1}) = 
\sum_{j=1}^{n-1} p_j(z_j)$ and $p_j$ are subharmonic, nonharmonic polynomials.
A first step in this direction is the $n=1$ case in $\C^2$. In
this case, $M$ is called a polynomial model and Nagel and Stein \cite{NaSt00} 
prove rapid decay of the heat kernel. The results in this article and \cite{Rai05h}
provide machinery to
improve the decay estimate to exponential decay. 
This is the topic of a forthcoming paper.

If $M$ is a polynomial model, then $M\cong \C\times\R$ and the realization of $\dbar_b$ (defined on $M$) on
$\C\times\R$ is $\bar L = \frac{\p}{\p\z} - i\frac{\p p}{\p\z}\frac{\p p}{\p t}$. $\bar L$ is translation
invariant in $t$, and taking a partial Fourier transform in $t$ produces $\Zbstp$. Thus, our understanding
of $e^{-s\Boxtp}$ for $\tau\in\R$ will yield results about $e^{-s\Box_b}$ on $M$. 
When $n>1$, and $M$ is a decoupled domain, the decoupling causes $\Box_{\tau P}$ to act diagonally.
If $\vartheta_q$ is the set of ordered $q$-tuples $(j_1,\dots,j_q)$ where
$1\leq j_1<j_2<\cdots<j_q\leq n$, it turns out that
$\Box_{\tau P} \big( \sum_{J\in \vartheta_q} \vp_J \, d\z_J\big) = \sum_{J\in\vartheta_q} \Box_J(\vp_J)\, d\z_J$,
and  each $\Box_J$ has a further decomposition into a sum of $(n-1)$ commuting operators of the form 
$\Boxtp$ or $\Boxwtp$. Thus, the study of $\Box_{\tau P}$ on $(0,q)$-forms reduces to the study
of the $\binom{n-1}q$ operators $\Box_J$ for $J\in\vartheta_q$. Furthermore, the commutation of the
elements $\Box_J$ makes the study of the heat operators $e^{-s\Boxtp}$ and $e^{-s\Boxwtp}$ critical
because if $\Box_J = \Box_1+\cdots+\Box_{n-1}$, then
$e^{-s\Box_J} = \prod_{j=1}^{n-1} e^{-s\Box_j}$.

\subsection{Refinement of the Problem and a Discussion of $\Boxwtp$}\label{subsec:refine}

When $\tau>0$, we \cite{Rai05h} write the solution
$u(s,z)$  of \eqref{eq:he i1}
as
\[
u(s,z) = \int_{\C} \Htp(s,z,w) f(w)\, dA(w),
\]
where  $\Htp(s,z,w)$ is a distributional kernel.
We show that  $\Htp$ is smooth (on an appropriate region)
and establish pointwise estimates for the kernels. We will extend our previous work and 
analyze the $\tau<0$ case. 

$\Boxtp$ has a special relationship with $\Boxwtp$.
\begin{align}
\Boxtp &= -\frac{\p^2 }{\p z\p\z} + \tau \frac{\p^2 p}{\p z\p\z}
+ \tau^2 \frac{\p p}{\p z}\frac{\p p}{\p\z}
+\tau\left( \frac{\p p}{\p z} \frac{\p }{\p\z} - \frac{\p p}{\p\z}\frac{\p }{\p z}\right) \label{eqn:Boxz}\\
&= -\frac14 \triangle + \frac14 \tau \triangle p+ \frac {\tau^2}4|\nabla p|^2  + \frac i2\tau 
\left(\frac{\p p}{\p x_1}\frac{\p}{\p x_2} - \frac{\p p}{\p x_2}\frac{\p}{\p x_1}\right) \label{eqn:Boxx}
\end{align}
while
\begin{align}
\Boxwtp &=\frac{\p^2 }{\p z\p\z} - \tau \frac{\p^2 p}{\p z\p\z}
+ \tau^2 \frac{\p p}{\p z}\frac{\p p}{\p\z}
+\tau\left( \frac{\p p}{\p z} \frac{\p }{\p\z} - \frac{\p p}{\p\z}\frac{\p }{\p z}\right) \label{eqn:Boxwz}\\
&= -\frac14 \triangle - \frac14 \tau \triangle p+ \frac {\tau^2}4|\nabla p|^2  + \frac i2\tau 
\left(\frac{\p p}{\p x_1}\frac{\p}{\p x_2} - \frac{\p p}{\p x_2}\frac{\p}{\p x_1}\right). \label{eqn:Boxwx}
\end{align}
The immediate consequence of \eqref{eqn:Boxx} and \eqref{eqn:Boxwx} is that if $\tilde p(x_1,x_2) = p(x_2,x_1)$, then
$\Box_{(-\tau)p} = \Boxw_{\tau \tilde p}$. Thus, our problem reduces from understand $\Boxtp$ and 
$\Boxwtp$ for all $\tau$ to understanding $\Boxtp$ and $\Boxwtp$ for $\tau>0$. Thus, the focus of this
paper is to understand (\ref{eq:he2}) when $\tau>0$.

$\Boxwtp$ is inherently more difficult to analyze than $\Boxtp$.
As a Schr\"odinger operator, $\Boxwtp$ has a 
nonpositive, unbounded potential. From the point of view of parabolic operator theory, the
$0^{\text{th}}$ order term can be negative and unbounded. In contrast,
$\Boxtp$ has a nonnegative potential and a nonnegative $0^{\text{th}}$ order term.
A further complication is that for $\tau>0$,
$\Null\Boxwtp\neq \{0\}$ (and in fact may be infinite dimensional, see \cite{Christ91}) 
while $\Null\Boxtp =\{0\}$. In fact, since $\Boxwtp$ has nonnegative eigenvalues
and is self-adjoint, it follows from the spectral theorem that
$\lim_{s\to\infty}e^{-s\Boxwtp} = \Ss$ where $\Ss$ is the Szeg\"o 
projection, i.e., the projection of $L^2(\C)$ onto $\Null\Zbstp$.
A consequence of the nonzero limit is that the kernel of $e^{-s\Boxwtp}$ cannot vanish
as $s\to\infty$. Thus, $\int_{0}^\infty e^{-s\Boxwtp}\,ds$ diverges and cannot be the
relative fundamental solution of $\Boxwtp$. $e^{-s\Boxwtp}(I-\Ss)$ functions as the
natural replacement for $e^{-s\Boxwtp}$ since $\int_{0}^\infty e^{-s\Boxwtp}(I-\Ss)\, ds$
does converge and equals the relative fundamental solution of $\Boxwtp$. Thus, we wish
to understand the pointwise size estimates of the kernel
of $e^{-s\Boxwtp}(I-\Ss)$ and its derivatives. 
Specifically, we will write
\[
e^{-s\Boxwtp}(I-\Ss)[f](z) = \int_{\C}\Gwtp(s,z,w)f(w)\, dA(w)
\]
and will analyize $\Gwtp(s,z,w)$ and its derivatives.

%
%
%
\section{Definitions and Principal Results}\label{sec:defns and results}
Let $p$ be a subharmonic, nonharmonic polynomial of degree $2m$. Let $D(z,r)$ be the 
Euclidean disk centered at $z$ of radius $r$.
Set
\begin{equation}\label{eq:Ajk}
\A{jk}z = \frac{1}{j! k!}\frac{\p^{j+k}p}{\p z^j\p\z^k}(z).
\end{equation}
Then $p(w) = \sum \A{jk}z(w-z)^j\overline{(w-z)}^k$. 
We define following  two  ``size" functions from the Carnot-Carath\'eodory geometry (see \cite{NaStWa85})
on the polynomial model $M_p = \{(z_1,z_2)\in\C^2 : \Imm z_2 = p(z_1)\}$. Since $\Zbstp$ is a partial
Fourier transform of $\bar L$, it not surprising that geometric elements on $M_p$ play a role in our analysis. 
Let
\begin{align}\label{eq:Lam}
\Lambda(z,\delta) &= \sum_{j,k\geq 1}  \left|\A{jk}z \right||\delta|^{j+k}
\intertext{and}
\label{eq:mu}
\mu_p(z,\delta) &= \inf_{j,k\geq 1} \Big|\frac{\delta}{\A{jk}z} \Big|^{\frac1{j+k}}.
\end{align}
The functions also arise
in the analysis of  magnetic Schr\"odinger  operators with electric potentials 
\cite{Shen96,Shen99,Kurata00}. 
It follows 
$\mu(z,\delta)$ is an approximate inverse to $\Lambda(z,\delta)$ if $\delta>0$. 
This means that if $\delta>0$,
\[
\mu\big(z,  \Lambda(z,\delta)\big) \sim \delta
\text{ and }
\Lambda\big(z,\mu(z,\delta)\big) \sim \delta.
\]
We use the notation $a\les b$ if $a\leq C b$ where $C$ is a constant that may depend on the
dimension 2 and the degree of $p$. We say that $a\sim b$ if $a\les b$ and $b\les a$.
For a little intuition for what $\mu_p(z,\delta)$ can look like, if
$x = \Rre z$ and $p(z) = x^{2m}$, then 
$\mu_p(z,\delta) \sim \frac{\delta^{1/2}}{x^{m-1}} + \delta^{\frac {1}{2m}}$.

We need to establish notation for adjoints. 
If $T$ is an operator (either 
bounded or closed and densely defined) on a Hilbert space with inner product
$\big(\,\cdot\, ,\cdot\,\big)$, let $T^*$ be the Hilbert space adjoint of $T$. This means that if
$f\in\Dom T$ and $g\in \Dom{T^*}$, then $\big(f,Tg\big) = \big(T^*f,g\big)$. 
If $U$ is an unbounded domain in some Euclidean space and
$T$ is an operator acting on $C^\infty_c(U)$ or $\mathcal{S}(U) = \{\vp\in C^\infty(U) : 
\vp \text{ has rapid decay}\}$, then we denote $T\sh$ as the adjoint in the sense of distributions.
This means that if $K$ is a distribution or a Schwartz distribution, then
$\langle T\sh K, \vp\rangle = \langle K, T \vp \rangle$. Note that if $T$ is not 
$\R$-valued, $T^* \neq T\sh$.

Since $\Boxtp$ is self-adjoint in $L^2(\C)$, it will follow that $\Htp(s,z,w) = \overline{\Htp(s,w,z)}$.
The relevance of this fact for adjoints 
is that the operators $\Zbstpw\sh$ and $\Zbstpw\sh$
are the appropriate differential operators to apply to $\Htp(s,z,w)$
and $\Zbstp\sh \neq\Zbstp$ or $\Zstp$. Thus, 
associated to a polynomial $p$ and the parameter $\tau\in\R$ are the weighted differential
operators 
\begin{align*}
\Wbstpw &= \frac{\p}{\p \w} -  \tau \frac{\p p}{\p \w}= e^{\tau p}\frac{\p p}{\p \w}e^{-\tau p}
& \Wstpw &= \frac{\p}{\p w} +  \tau \frac{\p p}{\p w}= e^{-\tau p}\frac{\p p}{\p w}e^{\tau p}.
\end{align*}
Observe that
$\overline{(\Zstp)} =  \Wbstp$ and $\overline{(\Zbstp)} = \Wstp$ and
\[
\Zbstp\sh = -\Wbstp \qquad\text{and}\qquad \Zstp\sh = -\Wstp.
\]

Analogously to writing $\frac{\p}{\p\z} = \frac 12\big(\frac{\p}{\p x_1} + i \frac{\p}{\p x_2}\big)$,
we let $X_1$ and $X_2$ denote the ``real" and ``imaginary" parts of $\Zstp$, that is,
\begin{align*}
X_1 &= \Zstp + \Zbstp = \frac{\p}{\p x_1} + i \tau \frac{\p p}{\p x_2} &
X_2 &= i(\Zstp - \Zbstp) = \frac{\p }{\p x_2} - i \tau \frac{\p p}{\p x_1}.
\end{align*}
Similarly, we let $U_1$ and $U_2$ denote the ``real" and ``imaginary" parts of $\Wstp$. Define
\begin{align*}
U_1 &= \Wstp + \Wbstp = \frac{\p}{\p x_1} - i \tau \frac{\p p}{\p x_2} &
U_2 &= i(\Wstp - \Wbstp) = \frac{\p }{\p x_2} +i \tau \frac{\p p}{\p x_1}.
\end{align*}
If $\alpha$ is a multiindex, we will also use the notation that $X^\alpha$ is a product of $|\alpha|$ 
operators of the form $X_1$ and $X_2$ while $U^\alpha$ is a similar product, except with $U_1$ and $U_2$
replacing $X_1$ and $X_2$.

We are now ready to state the main results. 

\begin{thm}\label{thm:H wig and deriv} Let $p$ be a subharmonic, nonharmonic polynomial and
$\tau>0$ a parameter. If
$n\geq 0$ and $Y^\alpha$ is a composition of $|\alpha|$
operators of the form $Y = \Zbstpz, \Zstpz$, $\Wbstpw, \Wstpw$, then there exist positive constants
$c, C_{|\alpha|}, C_{n,|\alpha|}$, so that
\[
\left|Y^\alpha \Hwtp(s,z,w)\right|
\leq C_{|\alpha|} e^{-c \frac{|z-w|^2}s} \max \left\{
\frac{e^{-c \frac{s}{\mu_p(w,1/\tau)^2}}e^{-c \frac{s}{\mu_p(z,1/\tau)^2}}} {s^{1+\frac 12|\alpha|}},
\frac{e^{-c \frac{|z-w|}{\mu_p(z,1/\tau)}}
e^{-c \frac{|z-w|}{\mu_p(w,1/\tau)}}}{\mu_p(w,1/\tau)^{2+|\alpha|}} \right\}
\]
Also, in the cases in which the derivatives 
annihilate the Szeg\"o kernel, i.e., $\frac{\p^n}{\p s^n}Y^\alpha \Ss(z,w)=0$,
the estimate simplifies to
\[
\left|\frac{\p^n}{\p s^n}Y^\alpha \Hwtp(s,z,w)\right|
\leq \frac{C_{n,|\alpha|}}{s^{1+n+\frac 12|\alpha|}} e^{-c \frac{|z-w|^2}s}
e^{-c \frac{s}{\mu_p(w,1/\tau)^2}}e^{-c \frac{s}{\mu_p(z,1/\tau)^2}}.
\]
\end{thm}

The corresponding theorem for estimates of $\Gwtp(s,z,w)$ is:
\begin{thm}\label{thm:G wig and deriv} Let $p$ be a subharmonic, nonharmonic polynomial and
$\tau>0$ a parameter. If
$n\geq 1$ and $Y^\alpha$ is a composition of
$|\alpha|$ operators of the form $Y= \Zbstpz, \Zstpz$, $\Wbstpw, \Wstpw$, then there exist 
positive constants
$c, C_{|\alpha|},  C_{|\alpha|,n}$ so that
\[
\left|Y^\alpha \Gwtp(s,z,w)\right|
\leq C_{|\alpha|} e^{-c \frac{s}{\mu_p(w,1/\tau)^2}}e^{-c \frac{s}{\mu_p(z,1/\tau)^2}} \max \left\{
\frac{ e^{-c \frac{|z-w|^2}s}} {s^{1+\frac 12|\alpha|}},
\frac{e^{-c \frac{|z-w|}{\mu_p(z,1/\tau)}}
e^{-c \frac{|z-w|}{\mu_p(w,1/\tau)}}}{\mu_p(w,1/\tau)^{2+|\alpha|}} \right\}
\]
Also, in the cases in which the derivatives annihilate the Szeg\"o kernel, i.e., 
$\frac{\p^n}{\p s^n}Y^\alpha \Ss(z,w)=0$, then $\frac{\p^n}{\p s^n}Y^\alpha \Hwtp(s,z,w)
= \frac{\p^n}{\p s^n}Y^\alpha \Gwtp(s,z,w)$, and
the estimate simplifies to
\[
\left|\frac{\p^n}{\p s^n}Y^\alpha \Hwtp(s,z,w)\right|
\leq \frac{C_{n,|\alpha|}}{s^{1+n+\frac 12|\alpha|}} e^{-c \frac{|z-w|^2}s}
e^{-c \frac{s}{\mu_p(w,1/\tau)^2}}e^{-c \frac{s}{\mu_p(z,1/\tau)^2}}.
\]
\end{thm}
Using techniques similar to those in \cite{NaSt00, Rai05h},
we have already
showed that $\Gwtp, \Hwtp \in C^\infty\big((0,\infty)\times\C\times\C)$ \cite{Rai05}. Our goal is to find
pointwise bounds on $\Hwtp$ and $\Gwtp$ and their derivatives.

\begin{rem} The Szeg\"o kernel $\Ss$ projects $L^2(\C)$ onto the kernel of $\Zbstp$, so
$\Zbstpz \Stp(z,w)=0$. Since $\Stp$ is self-adjoint in $L^2(\C)$, $\Wstpw \Stp(z,w)=0$ as well.
Moreover,  $\Stp(z,w)$ has no time dependence, so  $\frac{\p}{\p s}\Stp(z,w)=0$ which explains the lack
of $s$-derivatives in the first estimate of Theorem \ref{thm:H wig and deriv} and Theorem
\ref{thm:G wig and deriv}.
Also, using the techniques of \cite{Rai05h, NaSt00}, the self-adjointness
of $e^{-s\Boxwtp}$ and $\Ss$ implies
$\Hwtp(s,z,w) = \overline{\Hwtp(s,w,z)}$ and $\Gwtp(s,z,w) = \overline{\Gwtp(s,w,z)}$ 
(see \cite{Rai05} for details).
\end{rem}

To help us
understand Theorem \ref{thm:H wig and deriv} and Theorem \ref{thm:G wig and deriv}, let us recall 
results from \cite{Rai05h}.
\begin{thm}\label{thm:Rai05h est} 
Let $p$ be a subharmonic, nonharmonic polynomial and $\tau>0$ a parameter. If $n\geq 0$ and
$Y^\alpha$ is a product of $|\alpha|$ operators $Y = \Zbstpz,\Zstpz$ or $\Wbstpw, \Wstpw$, there exist constants
$C_{n,|\alpha|}$,  $c$ (independent of $\tau$) so that
\[
\left| \frac{\p^n}{\p s^n} Y^\alpha \Htp(s,z,w)\right| \leq \frac{C_{n,|\alpha|}}{s^{1+n+ \frac 12|\alpha|}}
e^{-c\frac{|z-w|^2}s} e^{-c \frac{s}{\mu_p(z, 1/\tau)^2}}e^{-c \frac{s}{\mu_p(w,1/\tau)^2}}.
\]
\end{thm}
By integrating in $s$, we have the corollary
\begin{cor}\label{cor:R,S decay} Let $\tau>0$.
Let $Y^\alpha$ be a product of $|\alpha|$ operators $Y = \Zbstpz,\Zstpz$ or $\Wbstpw, \Wstpw$. Then
there exist constants $c, C_{|\alpha|}>0$ so that
\begin{gather*}
|Y^\alpha \Rtp(z,w)| \leq C_{|\alpha|} 
 \begin{cases}
 |z-w|^{-1-|\alpha|} & |z-w|\leq \mu_p(z, \tfrac 1\tau)\\
 \frac{1}{\mu_p(z,1/\tau)^{1+|\alpha|}} e^{-c \frac{|z-w|}{\mu_p(z,1/\tau)}} 
e^{-c \frac{|z-w|}{\mu_p(w,1/\tau)}} 
                   &|z-w|\geq \mu_p(z,\tfrac 1\tau) 
 \end{cases}
\intertext{and}
|Y^\alpha \Ss(z,w)| \leq \frac{C_{|\alpha|}}{\mu_p(z,1/\tau)^{2+|\alpha|}}
e^{-c \frac{|z-w|}{\mu_p(z,1/\tau)}}e^{-c \frac{|z-w|}{\mu_p(w,1/\tau)}} .
\end{gather*}
\end{cor}

$e^{-s\Boxwtp}$ should behave well near $s=0$, and it does;  the estimates for $\Hwtp(s,z,w)$
agree with the estimates
for $\Htp(s,z,w)$ near 0. $e^{-s\Boxwtp}(I-\Ss)$ should behave well as $s\to\infty$, and we see that
$\Gwtp(s,z,w)$ has the 
same exponential decay in $s$ that $\Htp(s,z,w)$ has. On the other hand, 
$\lim_{s\to\infty} e^{-s\Boxwtp} = -\lim_{s\to 0}e^{-s\Boxwtp}(I-\Ss) = \Ss$, and the estimate
for $\Hwtp(s,z,w)$ as $s\to\infty$ and the estimate for $\Gwtp(s,z,w)$ near $s=0$ agree with the estimate
for $\Ss(z,w)$. Thus, the decay terms in Theorem \ref{thm:H wig and deriv} and 
Theorem \ref{thm:G wig and deriv}
are reasonable. In the case when the derivative annihilates $\Ss(z,w)$, we find that
the estimates for the derivatives of $\Hwtp(s,z,w)$ and $\Gwtp(s,z,w)$ agree and behave similarly to derivatives
of $\Htp(s,z,w)$.

%
%
%
%
\section{Embedding Theorems and Cancellation Conditions}\label{sec:embed and cancel}
In Section \ref{sec:Gw est}, we will derive estimates of $Y^\alpha \Gwtp(s,z,w)$ from results
about $\Htp(s,z,w)$. First, however, we need to prove some cancellation conditions for
$\Htp^s$ where $\Htp^s[\vp] = e^{-s\Boxtp}[\vp]$. 
In addition to their importance in this paper, the cancellation conditions complement
the size estimates of Theorem \ref{thm:Rai05h est} and Corollary \ref{cor:R,S decay}, and we feel
the results in this section are interesting in their own right.

We start by strengthening the Sobolev-embedding theorem from \cite{Rai05h}.
Define
\[
\tnabla = (X_1, X_2).
\]
Note $\tnabla = \nabla + i(\tau\frac{\p p}{\p x_2}, -\tau\frac{\p p}{\p x_1})$. To prove our 
weighted Sobolev embedding
theorem, we adapt ideas  from \cite{Ad78}.
\begin{thm}\label{thm:Poincare} Let $\Delta = (a_1,b_1)\times(a_2,b_2)\subset \R^2$ be a square
of sidelength $\delta$. If $(x_1,x_2)\in \Delta$ and if $f\in\mathcal{C}^2(\Delta)$, 
then
\[
|f(x_1,x_2)| ^2 \leq 4 \left( \frac{1}{\delta^2} \int_\Delta |f|^2 + \int_\Delta |\tnabla f|^2 
+ \delta^2 \int_\Delta |X_2X_1 f|^2\right).
\]
\end{thm}

\begin{proof} Our first goal is to show that for $\zeta\in (a_2,b_2)$,
\begin{equation}\label{eq:1d sobo}
\int_{a_1}^{b_1} |f(y_1,\zeta)|^2\, dy_1
 \leq 2 \left( \frac{1}{\delta}\int_\Delta |f|^2 + \delta \int_\Delta |X_2 f|^2\right).
\end{equation}
For a continuous function $u(y_1,y_2)$, by the Mean Value Theorem for Integrals, there exists
$\sigma\in (a_2,b_2)$ so that 
\begin{equation}\label{eqn:MVT for int}
\int_{\Delta} |u|^2 = \delta \int_{a_1}^{b_1} |u(y_1,\sigma)|^2\, dy_1.
\end{equation}
By the Fundamental Theorem of Calculus and Cauchy-Schwarz, for $\zeta \in (a_2,b_2)$
(without loss of generality $\zeta\geq \sigma$)
\begin{equation}\label{eqn:FTC u}
|u(y_1,\zeta)|^2 \leq 2\left( |u(y_1,\sigma)|^2 + \delta 
\int_\sigma^\zeta \left|\frac{\p u}{\p x_2}(y_1,t)\right|^2\, dt \right).
\end{equation}
We would like to replace the ordinary derivative in \eqref{eqn:FTC u} with $X_2$. To this end,
we define a function $A(y_1,y_2)$ so that $A(y_1,y_2) = -\int_\sigma^{y_2} \tau\frac{\p p}{\p y_1}(y_1,t)\,dt$.
Then $\frac{\p A}{\p y_2}(y_1,y_2) = -\tau\frac{\p p}{\p y_1}(y_1,y_2)$, and if we let
$u(y_1,y_2) = e^{iA(y_1,y_2)}f(y_1,y_2)$,
\[
\frac{\p}{\p y_2}\big( e^{iA(y_1,y_2)} f(y_1,y_2)\big) 
= e^{iA(y_1,y_2)} \left(\frac{\p f}{\p y_2} - i \tau\frac{\p p}{\p y_1} f\right)(y_1,y_2)
= e^{iA(y_1,y_2)}  X_2 f(y_1,y_2).
\]
Thus, \eqref{eqn:FTC u} can be rewritten as
\begin{equation}\label{eqn:FTC f}
|f(y_1,\zeta)|^2 \leq 2\left( |f(y_1,\sigma)|^2 + \delta\int_\sigma^\zeta |X_2 f(y_1,t)|^2\, dt\right)
\leq 2\left( |f(y_1,\sigma)|^2 + \delta\int_{a_2}^{b_2} |X_2 f(y_1,t)|^2\, dt\right).
\end{equation}
Integrating \eqref{eqn:FTC f} over $(a_1,b_1)$ and using \eqref{eqn:MVT for int}, we have
\begin{align*}
\int_{a_1}^{b_1} |f(y_1,\zeta)|^2\, dy_1 
&\leq 2\left( \int_{a_1}^{b_1} |f(y_1,\sigma)|^2\, dy_1 + \delta\int_{\Delta} |X_2 f|^2 \right)\\
&\leq 2\left( \frac{1}{\delta}\int_\Delta |f|^2  + \delta\int_{\Delta} |X_2 f|^2\right)
\end{align*}
which is \eqref{eq:1d sobo}.

We now finish the proof. Given $u\in \mathcal{C}^2(\Delta)$,
by the Mean Value Theorem for Integrals, there exists $\varpi\in (a_1,b_1)$
so that
\begin{equation}\label{eqn:MVT 2}
\int_{a_1}^{b_1} |u(s,x_2)|^2\, ds = \delta |u(\varpi,x_2)|^2.
\end{equation}
By the Fundamental Theorem of Calculus and H\"older's inequality (without loss of generality
$x_1>\varpi$),
\begin{equation}\label{eqn:FTC 2}
|u(x_1,x_2)|^2 \leq 2\left( |u(\varpi,x_2)|^2 
+ \delta \int_{\varpi}^{x_1} \left|\frac{\p u}{\p x_1}(s,x_2)\right|^2\, ds\right).
\end{equation}
Similarly to above, we let $B(x_1,x_2) = \int_{\varpi}^{x_1} \tau\frac{\p p}{\p x_2}(s,x_2)\, ds$,
so $\frac{\p B}{\p x_1}(x_1,x_2) = \tau\frac{\p p}{\p x_2}(x_1,x_2)$. If we set
$u(x_1,x_2) = e^{iB(x_1,x_2)}f(x_1,x_2)$, then analogously to earlier,
\[
\frac{\p u}{\p x_1}(x_1,x_2) = e^{iB(x_1,x_2)} X_1 f(x_1,x_2).
\]
Thus, substituting $e^{iB(x_1,x_2)}  f(x_1,x_2)$ into \eqref{eqn:FTC 2} and using
\eqref{eqn:MVT 2} and \eqref{eq:1d sobo}, we have
\begin{align*}
|f(x_1,x_2)|^2 &\leq 2\big( |f(\varpi,x_2)|^2 + \delta\int_{a_1}^{b_1}|X_1 f(s,x_2)|^2\, ds \big) \\
&= 2\left( \frac 1\delta \int_{a_1}^{b_1} |f(s,x_2)|^2\, ds 
+ \delta\int_{a_1}^{b_1}|X_1 f(s,x_2)|^2\, ds \right) \\
&\leq 4\left( \frac{1}{\delta^2} \int_\Delta |f|^2 + \int_\Delta |X_2 f|^2 + \int_\Delta |X_1 f|^2
+ \delta^2 \int_\Delta |X_2X_1 f|^2\right).
\end{align*}
\end{proof}

From \cite{Rai05h}, we can write $X_j X_k = A_{jk}\Boxtp$ where $A_{jk}$ is a bounded operator in
$L^2(\C)$, and we can estimate 
$\|X f\|_{L^2(\C)} \les \frac 1\delta \|f\|_{L^2(\C)} + \delta \|\Boxtp f\|_{L^2(\C)}$.
Combining these facts with Proposition \ref{thm:Poincare} and \cite{Rai05h}, we have the  corollary:
\begin{cor}\label{cor:sobolev} Let $\tau>0$.
There exists $C>0$ so that if $f\in \mathcal{C}^2(D(z,\delta))$ and $\supp f \in D(z,\delta)$, then
\[
|f(z)| \leq \frac C\delta \left( \|f\|_{L^2(D(z,\delta))} + \delta^2 \|\Boxtp f\|_{L^2(D(z,\delta))}\right).
\]
\end{cor}

We need the following Poincar\'e lemma.
\begin{prop} \label{prop:L2 vs. L2 deriv}
Let $\vp\in \mathcal{C}^1_{c}(\Delta)$ where
$\Delta = (a_1,b_1)\times(a_2,b_2)\subset\R^2$ is a square with sidelength $\delta$.  Then
\[
\|\vp\|_{L^2(\Delta)} \leq \sqrt 2 \delta \|X_j\vp \|_{L^2(\C)}
\]
where $j=1$ or $2$.
\end{prop}

\begin{proof}By the Fundamental Theorem of Calculus and Cauchy-Schwarz,
\[
|f(x_1,x_2)|^2 \leq 2\Big(|f(x_1,\zeta)|^2 + \delta \int_{a_2}^{b_2} \Big|\frac{\p f}{\p x_2}(x_1,t)\Big|^2\,dt\Big).
\]
Choosing $\zeta$ outside of $(a_2,b_2)$ and integrating over $\Delta$ (in $x_1$ and $x_2$), we see
that
\begin{equation}\label{eqn:f L2 vs. Df L2}
\|f\|_{L^2(\Delta)}^2 \leq 2\delta^2 \Big \|\frac{\p f}{\p x_2}\Big\|_{L^2(\Delta)}^2.
\end{equation}
As in the proof of Theorem \ref{thm:Poincare}, define a function
$A(x_1,t) = -\int_{0}^{t} \tau\frac{\p p}{\p x_1}(x_1,s)\, ds$ so that $\frac{\p A}{\p t}(x_1,t) = 
-\tau\frac{\p p}{\p x_1}(x_1,t)$. If we let $f(x_1,x_2) = e^{iA(x_1,x_2)} \vp(x_1,x_2)$, then
\begin{equation}\label{eqn:f to vp}
\Big|\frac{\p f}{\p x_2}\Big| = \Big|\frac{\p \vp}{\p x_2} - i\tau\frac{\p p}{\p x_1}\Big| = \big|X_2 \vp\big|, 
\end{equation}
and plugging in (\ref{eqn:f to vp}) into (\ref{eqn:f L2 vs. Df L2}) finishes the proof for $j=2$. The $j=1$ case
is analogous and proven by integrating in $x_1$ first.
\end{proof}

The next lemma proves  a cancellation condition for the heat kernel $\Htp(s,z,w)$. 
\begin{lem}\label{lem:H cancel}
Let $\tau>0$ and $\delta > 0$ and
$\vp\in \cic{D(z,\delta)}$. If $Y^\alpha$ is a composition of $|\alpha|$ operators
of the form $Y = \Zstp$ and $\Zbstp$, then there exists a constant $C_{|\alpha|}$ so that if $|\alpha|=2k$ is
even,
\[
|Y^\alpha \Htp^s[\vp](z)| 
\leq C_{|\alpha|} \delta^{-1}
\big(\|\Boxtp^k[\vp]\|_{L^2(\C)} + \delta^{2} \|\Boxtp^{k+1}[\vp]\|_{L^2(\C)}\big)
\]
and if $|\alpha|=2k+1$ is odd,
\[ 
|Y^\alpha \Htp^s[\vp](z)| 
\leq C_{|\alpha|} \delta^{-1}\big(\delta\|\Boxtp^{k+1}\vp\|_{L^2(\C)} + \delta^{3}\|\Boxtp^{k+2}\vp\|_{L^2(\C)}\big)
\]
\end{lem}

\begin{proof}
From \cite{Rai05h} and \cite{Rai05f}, if $|\alpha|=2k$ is even, then
$\| Y^\alpha \vp \|_{L^2(\C)} \leq C_{\gamma,p} \|\Boxtp^k \vp \|_{L^2(\C)}$.
By Corollary \ref{cor:sobolev},  it follows that
\begin{align*}
|Y^\alpha \Htp^s[\vp](z)| 
&\leq \frac{C}{\delta}\big(\|Y^\alpha\Htp^s[\vp]\|_{L^2(\C)} + \delta^2 \|\Boxtp Y^\alpha \Htp^s[\vp]\|_{L^2(\C)}\big)\\
&\leq \frac{C}{\delta}\big(\|\Htp^s\Boxtp^k[\vp]\|_{L^2(\C)} + \delta^2 
\|\Htp^s\Boxtp^{k+1}[\vp]\|_{L^2(\C)}\big) \\
&\leq \frac{C}{\delta}\big(\|\Boxtp^k[\vp]\|_{L^2(\C)} + \delta^2 
\|\Boxtp^{k+1}[\vp]\|_{L^2(\C)}\big).
\end{align*}
The last line uses the fact that $e^{-s\Boxtp}$ is a contraction.
The $|\alpha|=2k+1$ case is similar. If $|\beta|= 2n+1$ for some nonnegative integer $n$, then
\begin{align*}
\| X^\beta \Htp^s[\vp]\|_{L^2(\C)}
\leq C \sum_{j=1}^2 \|X_j\Htp^s\Boxtp^n[\vp]\|_{L^2(\C)} 
&\leq \frac C\delta \big(\|\Htp^s\Boxtp^n[\vp]\|_{L^2(\C)} + \delta^2 \|\Htp^s\Boxtp^{n+1}[\vp]\|_{L^2(\C)}\big)\\
&\leq \frac C\delta \big(\|\Boxtp^n[\vp]\|_{L^2(\C)} + \delta^2 \|\Boxtp^{n+1}[\vp]\|_{L^2(\C)}\big).
\end{align*}
Thus,
\begin{align*}
|Y^\alpha\Htp^s[\vp](z)|
&\leq \frac{C}{\delta}\big( \|Y^\alpha\Htp^s[\vp]\|_{L^2(\C)} + \|\Boxtp Y^\alpha \Htp^s[\vp]\|_{L^2(\C)}\big)\\
&\leq \frac{C}{\delta} \big(\delta^{-1}\|\Boxtp^k[\vp]\|_{L^2(\C)}
+\delta\|\Boxtp^{k+1}[\vp]\|_{L^2(\C)} + \delta^{3}\|\Boxtp^{k+2}[\vp]\|_{L^2(\C)}\big).
\end{align*}
Using Proposition \ref{prop:L2 vs. L2 deriv}, we are done.
\end{proof}

\begin{rem} It is possible to prove a slightly weaker version of
Lemma \ref{lem:H cancel} using the correspondence between NIS operators  
and the one-paramater families (OPF) of operators in \cite{Rai05f}.
Nagel and Stein
\cite{NaSt00} prove that the heat operators $e^{-s\Box_b}$ and $e^{-s\Box_b}(I-S)$ are 
nonisotropic smoothing (NIS)  operators
of order 0 with NIS constants independent of $s$.  Thus $\Htp^s$ is an
OPF operator of order 0 with OPF constants independent of $s$. A result would then follow immediately.
The problem, however, is that
we lose control over the number of derivatives needed to bound
$|Y^\alpha \Htp^s[\vp]|$. Also,
the definition of NIS operator in \cite{NaSt00} is different than the definition of NIS operator
in \cite{NaRoStWa89} on which the correspondence is based. Since the kernels of 
$e^{-s\Box_b}$ and $e^{-s\Box_b}(I-S)$ are smooth, however, it is trivial to show that 
$e^{-s\Box_b}$ and $e^{-s\Box_b}(I-S)$ are NIS operators in the sense of \cite{NaRoStWa89}. 
Regardless, Lemma \ref{lem:H cancel} is a direct and more natural proof.
\end{rem}

Recall that $G_{\tau p}$ is the solving operator for $\Boxtp$ when $\tau>0$.
\begin{lem}\label{lem:R cancel}
Let $\tau>0$.
Let $\delta>0$ and $\vp\in \cic{D(z,\delta)}$. If $Y^\alpha$ is a composition of $|\alpha|$ operators of the form
$Y = \Zbstp$ or $\Zbstp$, then there exists a constant $C_{|\alpha|}$ so that if $|\alpha|=2k>0$ is even
or $|\alpha|=0$ and $\delta \geq \mu_p(z,\frac 1\tau)$,
then
\[
|Y^\alpha G_{\tau p}[\vp](z)| 
\leq C_{|\alpha|} \delta \big(\|\Boxtp^k\vp\|_{L^2(\C)} +\delta^2 \|\Boxtp^{k+1}\vp\|_{L^2(\C)}\big)
\]
and if $|\alpha|=2k+1>0$ is odd, then
\[
|Y^\alpha G_{\tau p}[\vp](z)| 
\leq C_{|\alpha|} \delta
\big( \delta \|\Boxtp^{k+1}\vp\|_{L^2(\C)} + \delta^3 \|\Boxtp^{k+2}\vp\|_{L^2(\C)}\big).
\]
If $|\alpha|=0$ and $\delta < \mu_p(z, \frac 1\tau)$
\[
|G_{\tau p}[\vp](z)| \leq C_0 \delta \Big( \log(\tfrac{2\mu_p(z,\frac 1\tau)}{\delta}) \|\vp\|_{L^2(\C)}
+ \delta^2 \|\Boxtp\vp\|_{L^2(\C)}\Big).
\]
\end{lem}

\begin{proof}
Recall that $G_{\tau p}[\vp](z) = \int_{0}^\infty e^{-s\Boxtp}[\vp](z)\, ds$. Thus,
\begin{equation}
|Y^\alpha G_{\tau p}[\vp](z)| \leq \left|\int_0^{\delta^2} Y^\alpha e^{-s\Boxtp}[\vp](z)\, ds \right|
+ \left|\int_{\delta^2}^\infty Y^\alpha e^{-s\Boxtp}[\vp](z)\, ds \right|\label{eqn:R decomp}
\end{equation}
We can estimate the first integral in \eqref{eqn:R decomp} with Lemma \ref{lem:H cancel}. If
$|\alpha|=2k$, then 
\[
\left|\int_0^{\delta^2} Y^\alpha e^{-s\Boxtp}[\vp](z)\, ds \right|
\leq \delta^2 \frac {C_{|\alpha|}}{\delta}\big( \|\Boxtp^k\vp\|_{L^2(\C)} 
+ \delta^2 \|\Boxtp^{k+1}\vp\|_{L^2(\C)}\big).
\]
The first integral in \eqref{eqn:R decomp} for the 
$|\alpha|=2k+1$ case is handled similarly. 

For the tail, we first assume $|\alpha|\geq 1$.
\begin{align*}
\Big|\int_{\delta^2}^\infty &Y^\alpha e^{-s\Boxtp}[\vp](z)\, ds \Big|
\leq \int_{\delta^2}^\infty\int_\C |Y^\alpha_z \Htp(s,z,w)\vp(w)|\, dA(w) ds\\
&\leq \int_{\C}|\vp(w)|\int_{\delta^2}^\infty \frac{1}{s^{1+\frac 12|\alpha|}}\, ds
\les \delta\|\vp\|_{L^2(\C)}\delta^{-|\alpha|} = \delta^{1-|\alpha|}\|\vp\|_{L^2(\C)}.
\end{align*}
A repeated use of Proposition \ref{prop:L2 vs. L2 deriv} shows 
\[
\delta^{-|\alpha|}\|\vp\|_{L^2(\C)} \leq 2^{|\alpha|/2}\|X_2^{|\alpha|} \vp\|_{L^2(\C)},
\]
and results from \cite{Rai05h} allow to change $X_2^\alpha$ to powers of $\Boxtp$ (as done above).
If $|\alpha|=0$, then 
\begin{equation}\label{eqn:log G cancel est}
\Big|\int_{\delta^2}^\infty e^{-s\Boxtp}[\vp](z)\, ds \Big|
\les \int_{\delta^2}^\infty \int_\C \frac 1s e^{-c \frac{s}{\mu_p(z,1/\tau)^2}}\, ds 
\les \delta \|\vp\|_{L^2(\C)} \int_{\frac{\delta^2}{\mu_p(z,1/\tau)^2}}^\infty \frac 1r e^{-r}\, dr
\end{equation}
If $\delta\geq \mu_p(z,\frac 1\tau)$, then (\ref{eqn:log G cancel est}) is bounded by
$\delta \|\vp\|_{L^2(\C)}$, but if $\delta < \mu_p(z,\frac 1\tau)$, then estimating the final term
in (\ref{eqn:log G cancel est}),
\[
\int_{\frac{\delta^2}{\mu_p(z,\frac 1\tau)^2}}^\infty \frac 1r e^{-r}\, dr
 = \int_{\frac{\delta^2}{\mu_p(z,\frac 1\tau)^2}}^1 \frac 1r e^{-r}\, dr 
+ \int_{1}^\infty \frac 1r e^{-r}\, dr \leq 2\log(\tfrac{2\mu_p(z,1/\tau)}{\delta}).
\]
\end{proof}

%
%
\section{Estimating $\Gwtp(s,z,w)$}\label{sec:Gw est}
The key to estimating $\Gwtp(s,z,w)$ and its derivatives lies in the elementary statement that
\[
\Boxtp \Zbstp = -\Zbstp\Zstp\Zbstp = \Zbstp\Boxwtp.
\]
From this statement, we have its corollary, that
\[
\Boxtp^n\Zbstp = \Zbstp\Boxwtp^n
\]
and
\[
\frac{(-1)^n s^n \Boxtp^n}{n!}\Zbstp = \Zbstp \frac{(-1)^n s^n\Boxwtp^n}{n!}.
\]
By the spectral theorem, summing in $n$ proves the first statement of
\begin{prop}\label{prop:box boxw relationship}
\[
e^{-s\Boxtp}\Zbstp = \Zbstp e^{-s\Boxwtp}
\]
and
\[
e^{-s\Boxwtp}\Zstp = \Zstp e^{-s\Boxtp}.
\]
\end{prop}
\begin{proof}
The second statement follows from the equality $\Boxwtp\Zstp = -\Zstp\Zbstp\Zstp = \Zstp\Boxtp$.
\end{proof}

Translating Proposition \ref{prop:box boxw relationship} to the kernel side, we see that:

\begin{cor}\label{cor:integral relating Gwtp and Htp}
Let $\tau>0$. Let $n\geq 0$ be an integer and 
$Y_z^\alpha$ be a product of $|\alpha|$ operators of the form $Y = \Zbstpz, \Zstpz$
or $\Wbstpw, \Wstpw$. Then
\[
\frac{\p^n}{\p s^n}Y^\alpha  \Zbstpz\Hwtp(s,z,w) 
= \frac{\p^n}{\p s^n}Y^\alpha \Zbstpz\Gwtp(s,z,w) 
=\frac{\p^n}{\p s^n}Y^\alpha \Wbstpw\Htp(s,z,w), 
\]
and
\[
\frac{\p^n}{\p s^n}Y^\alpha  \Wstpw\Hwtp(s,z,w) 
= \frac{\p^n}{\p s^n}Y^\alpha  \Wstpw\Gwtp(s,z,w) 
=\frac{\p^n}{\p s^n}Y^\alpha  \Zstpz\Htp(s,z,w). 
\]
Also, if $Y^\alpha = X^\beta_z U^\gamma_w$, then
\[
X_z^\beta U^\gamma_w \Gwtp(s,z,w) 
= -\int_\C U^\gamma_w\Wbstpw\Htp(s,v,w) X^\beta_z \Rtp(z,v)\, dA(v).
\]
\end{cor}

\begin{proof}
The first two equalities follow immediately from interpreting Proposition \ref{prop:box boxw relationship}
on the kernel side. The third equality requires a short argument. Since
\[
\Zbstp e^{-s\Boxwtp} = e^{-s\Boxtp}\Zbstp,
\]
we can apply $\Rtp$ to both sides and get
\[
(I-\Ss)e^{-s\Boxwtp} = \Rtp e^{-s\Boxtp}\Zbstp.
\]
The justification that $\Rtp\Zbstp = I-\Ss$ is \cite{Christ91}, pp.219-220. 
But $(I-\Ss)e^{-s\Boxwtp} = \Gwtp^s$, and
\begin{align*}
\Rtp e^{-s\Boxtp}[\Zbstp\vp](z)
&= -\Rtp\left[\int_\C \Wbstpw\Htp(s,\cdot,w)\vp(w)\, dA(w)\right](z)\\
&= -\int_C\left(\int_\C \Rtp(z,v)\Wbstpw \Htp(s,v,w)\, dA(v)\right) \vp(w)\, dA(w).
\end{align*}
Thus
\[
\Gwtp(s,z,w) = -\int_\C \Wbstpw\Htp(s,v,w) \Rtp(z,v)\, dA(v).
\]
Applying derivatives finishes the proof.
\end{proof}

Corollary \ref{cor:integral relating Gwtp and Htp} and Theorem \ref{thm:Rai05h est} give immediate estimates
for special derivatives of $\Gwtp(s,z,w)$ and $\Hwtp(s,z,w)$.

%
%
\begin{thm}\label{thm:Gw Hw special derivatives}
Let $p$ be a subharmonic, nonharmonic polynomial and $\tau>0$ a parameter. Let
$n\geq 0$ be an integer and 
$Y^\alpha$ a product of $|\alpha|$ operators of the form $Y = \Zbstpz,\Zstpz$, or 
$\Wbstpw,\Wstpw$. If $\frac{\p^n}{\p s^n}Y^\alpha$ annihilates the 
Szeg\"o kernel, i.e., $\frac{\p^n}{\p s^n}Y^\alpha\Stp(z,w)=0$,  then there exist
constants $c, C_{n,|\alpha|}>0$ (depending on $n$ and $|\alpha|$) so that
\[
\left| \frac{\p^n}{\p s^n}Y^\alpha \Hwtp(s,z,w) \right| 
= \left|\frac{\p^n}{\p s^n}Y^\alpha\Gwtp(s,z,w)\right| 
\leq \frac{C_{|\alpha|}}{s^{1+n+\frac 12|\alpha|}} 
e^{-c \frac{|z-w|^2}s} e^{-c \frac{s}{\mu_p(w,1/\tau)^2}}e^{-c \frac{s}{\mu_p(z,1/\tau)^2}}.
\]
\end{thm}

The last result we need to prove Theorem \ref{thm:G wig and deriv} is  from \cite{Christ91} relating the size
of $\mu_p(z,\frac 1\tau)$ and $\mu_p(w,\frac 1\tau)$.
\begin{lem}\label{lem:3.3 Christ} If $\tau>0$ and $|z-w|>\mu_p(w,\frac 1\tau)$, then
\[
\frac{\mu_p(z,1/\tau)}{\mu_p(w,1/\tau)} +  
\frac{\mu_p(w,1/\tau)}{\mu_p(z,1/\tau)} 
\leq C \left(\frac{|z-w|}{\mu_p(w,1/\tau)}\right)^{\deg p}
\]
for some constant $C<\infty$.
\end{lem}

One effect of Lemma \ref{lem:3.3 Christ} is 
that in a situation where $|z-w|\geq \mu_p(z,\frac 1\tau)$ (which
implies that $|z-w|\ges \mu_p(w,\frac 1 \tau)$) we can perform the following estimate
\begin{equation}\label{eqn:mu(z) to mu(w)}
\frac{1}{\mu_p(z,\tfrac 1\tau)} e^{-c \frac{|z-w|}{\mu_p(z,1/\tau)}}
=\frac{1}{\mu_p(w,\tfrac 1\tau)} \frac{\mu_p(w,1/\tau)}{\mu_p(z,1/\tau)}
e^{-c \frac{|z-w|}{\mu_p(z,1/\tau)}}
\les \frac{1}{\mu_p(w,\tfrac 1\tau)} e^{-c \frac{|z-w|}{\mu_p(z,1/\tau)}}
\end{equation}
with a possible decrease of $c$. (\ref{eqn:mu(z) to mu(w)}) reiterates that the statements of
Theorem \ref{thm:H wig and deriv} and Theorem \ref{thm:G wig and deriv} are symmetric in 
$z$ and $w$.

\begin{proof}[Proof (Theorem \ref{thm:G wig and deriv}).]
We estimate $X^\alpha_z U^\beta_w\Gwtp(s,z,w)$. 
This estimate and
Theorem \ref{thm:Gw Hw special derivatives} prove Theorem  \ref{thm:G wig and deriv}.
We are now ready to estimate the integral in Corollary \ref{cor:integral relating Gwtp and Htp}.

It is enough to prove 
\[
|X^\alpha_z U^\beta_w\Gwtp(s,z,w)|
\les e^{-c \frac{s}{\mu_p(w,1/\tau)^2}} \max\bigg\{
\frac{e^{-c \frac{|z-w|^2}s}} {s^{1+\frac 12(|\alpha|+|\beta|)}}, 
\frac{e^{-c \frac{|z-w|}{\mu_p(z,1/\tau)}}}{\mu_p(w,1/\tau)^2} \bigg\} 
\]
because $\Gwtp(s,z,w)$ is conjugate symmetric in $z$ and $w$, and all of the estimates we use
in the proof of Theorem \ref{thm:G wig and deriv} are symmetric in $z$ and $w$.

Let $U^\gamma_w = U^\beta_w \Wbstpw$. Then $|\gamma|\geq 1$. For $x\in\C$, let $\vp_x\in \cic\C$ so that
$\vp_x \equiv 1$ on $D(x,\mu_p(x,\frac 1\tau))$, $\vp_x \equiv 0$ off 
$D(x,2\mu_p(x,1))$, and $|D^n\vp_x| \leq
c_n \mu_p(x,\frac 1\tau)^{-n}$.
For $|z-w|>2\max\{\mu_p(z,\frac 1\tau),\mu_p(w,\frac 1\tau)\}$, we will estimate
\begin{align*}
\int_\C U^\beta_w \Wbstpw\Htp(s,v,w) & X^\alpha_z \Rtp(z,v)\, dA(v)
= \hspace{-1.785pt}
\underbrace{\int_\C U^\gamma_w\Htp(s,v,w) X^\alpha_z\Rtp(z,v)\vp_z(v)(1-\vp_w(v))\, dA(v)}_I \\
&+ \underbrace{\int_\C U^\gamma_w\Htp(s,v,w) X^\alpha_z\Rtp(z,v)(1-\vp_z(v))\vp_w(v)\, dA(v)}_{II}\\
&+ \underbrace{\int_{|v-w|>|v-z|} U^\gamma_w\Htp(s,v,w) X^\alpha_z\Rtp(z,v)(1-\vp_z(v))
(1-\vp_w(v))\, dA(v)}_{III}\\
&+ \underbrace{\int_{|v-w|<|v-z|} U^\gamma_w\Htp(s,v,w) X^\alpha_z\Rtp(z,v)(1-\vp_z(v))
(1-\vp_w(v))\, dA(v).}_{IV}
\end{align*}

We start by estimating $|I|$. We use Lemma \ref{lem:R cancel}. Let $\Htp^{s,w}(v) = \Htp(s,v,w)$. 
If $|\alpha|=2k-1$ is odd and $\delta = \mu_p(z,\frac 1\tau)$, then
\begin{align*}
|I| &= | X^\alpha_z \Rtp[U^\gamma_w\Htp^{s,w} \vp_z](w)| \\
&\leq c_{|\alpha|}\Big(\delta\|\Boxtp^k\big(U^\gamma_w\Htp(s,v,w)\vp_z(v)\big)\|_{L^2(\C)}
+ \delta^3\|\Boxtp^{k+1}\big(U^\gamma_w\Htp(s,v,w)\vp_z(v)\big)\|_{L^2(\C)}\Big).
\end{align*}
The two terms are estimated similarly. We present an estimate of the first term.
\begin{align*}
&\delta\|\Boxtp^k\big(U^\gamma_w\Htp(s,v,w)\vp_z(v)\big)\|_{L^2(\C)}
\leq \delta \sum_{|I|\leq 2k} \|X_v^I U_w^\gamma\Htp(s,v,w) D^{2k-|I|}\vp_z(v)\|_{L^2(\C)}\\
&\les \delta \sum_{|I|\leq 2k} \frac{1}{s^{1+\frac 12|I|+\frac 12|\gamma|}}
\delta e^{-c\frac{|z-w|^2}s} e^{-c \frac{s}{\mu_p(w,1/\tau)^2}} \delta^{-2k+|I|} \delta \\
&\les \frac{1}{s^{1+\frac 12(|\alpha|+|\beta|)}} e^{-c\frac{|z-w|^2}s} e^{-c \frac{s}{\mu_p(w,1/\tau)^2}}.
\end{align*}
The last line follows immediately from the previous line if $s<\mu_p(w,\frac 1\tau)^2$. If
$s>\mu_p(w,\frac 1\tau)^2$, we can use the decay of $e^{-c \frac{s}{\mu_p(w,1/\tau)^2}}$ to achieve the
desired bound. Also, the $|\alpha|$ even case uses the same argument as the $|\alpha|$ odd case.

To estimate $II$, 
note that on the support of $\vp_w$, $|v-z|\sim|w-z|$ and $\mu_p(v,\frac 1\tau)
\sim \mu_p(w,\frac 1\tau)$. Assume first that $|\gamma|=2k$ is even.
If $s<\mu_p(w,\frac 1\tau)^2$ and we set 
$\Rtp^z(v) = \Rtp(z,v)$ and $\delta = \mu_p(w,\frac 1\tau)$, then by Lemma \ref{lem:H cancel}, 
\begin{align*}
|II| &= \left | U^\gamma_w (\Htp^s)\sh \big[ X^\alpha_z \Rtp^z \vp\big](w) \right| \\
&\leq \frac C\delta \Big( \big\|(\Box_{\tau p, v}\sh)^k\big(X^\alpha_z \Rtp^z \vp_w \big) \big\|_{L^2(\C)}
+ \delta^2 \big\|(\Box_{\tau p, v}\sh)^{k+1}\big(X^\alpha_z \Rtp^z \vp_w \big) \big\|_{L^2(\C)}\Big).
\end{align*}
We estimate the first term.
\begin{align*}
\frac 1\delta \big\|(\Box_{\tau p, v}\sh)^k\big(X^\alpha_z \Rtp^z \vp_w \big) \big\|_{L^2(\C)}
&= \frac 1\delta\sum_{|I|\leq 2k} \big\| (U^I_v X^\alpha_z \Rtp(z,v)) D^{2k-I}\vp_w(v)\big\|_{L^2(\C)} \\
&\les \frac 1\delta e^{-c\frac{|z-w|}{\mu_p(z,1/\tau)}} \frac{1}{\mu_p(z,1/\tau)^{|\alpha|+1}
\mu_p(w,\frac 1\tau)^{|I|}} \mu_p(w,1/\tau)^{|I|-|\gamma|} \delta \\
&\leq e^{-c\frac{|z-w|}{\mu_p(z,1/\tau)}} \frac{1}{\mu_p(z,1/\tau)^{|\alpha|+1}
\mu_p(w,1/\tau)^{|\beta|+1}}.
\end{align*}
The second term is handled similarly, as is the $|\gamma|$ odd case.
If $s \geq \mu_p(w,\frac 1\tau)^2$, we would like the estimate of $II$ as above except with an additional
$e^{-c\frac{s}{\mu_p(w,1/\tau)^2}}$ term. For this estimate, we just use the size conditions for
$\Htp$ and $\Rtp$.
\[
|II| \les e^{-c \frac{|z-w|}{\mu_p(z,1/\tau)}} \frac{1}{\mu_p(z,\frac 1\tau)^{1+|\alpha|}}
\frac{1}{s^{1+\frac12|\gamma|}}e^{-c\frac{s}{\mu_p(w,1/\tau)^2}} \mu_p(w,\tfrac 1\tau)^2
\leq \frac {e^{-c\frac{|z-w|}{\mu_p(z,1/\tau)}}e^{-c\frac s{\mu_p(w,1/\tau)^2}}} 
{\mu_p(z,1/\tau)^{|\alpha|+1} \mu_p(w,1/\tau)^{|\beta|+1}}.
\]

Estimating $III$ and $IV$ follows from the size estimates for $\Rtp$ and $\Htp$
from Theorem \ref{thm:Rai05h est} and Corollary \ref{cor:R,S decay}. Note
that if $|v-w|>|v-z|$, then $|v-w| \geq \frac 12|z-w|$ and
\begin{align}
|III|
&\leq \int_{\atopp{|v-w|>\mu_p(w,1/\tau)}{\atopp{|v-z|>\mu_p(z,1/\tau)}{|v-w|>|v-z|}}}
|U^\gamma_w\Htp(s,v,w)  X^\alpha_z\Rtp(z,v)|\, dA(v) \nn\\
&\les \frac{1}{s^{1+\frac 12|\gamma|}} e^{-c\frac{|z-w|^2}s} e^{-c \frac{s}{\mu_p(w,1/\tau)^2}}
\int_{|v-z|>\mu_p(z,1/\tau)} \frac{1}{\mu_p(z,\frac1\tau)^{|\alpha|+1}}
e^{-c \frac{|z-v|}{\mu_p(z,1/\tau)}}\, dA(v) \nn\\
&\leq \frac 1{s^{1 + \frac 12|\beta| + \frac 12}} \frac 1{\mu_p(z,1/\tau)^{|\alpha|-1}}
e^{-c\frac{|z-w|^2}s} e^{-c \frac{s}{\mu_p(w,1/\tau)^2}}.
\label{eqn:III est}
\end{align}
(\ref{eqn:III est}) is close the desired estimated, but we need to replace
$\frac{1}{\mu_p(z,1/\tau)^{|\alpha|-1}}$ with $\frac{1}{s^{\frac 12 |\alpha|-\frac 12}}$. There are two cases.
If $s>\mu_p(z,\frac 1\tau)^2$, then (shrinking $c$ as necessary) by Lemma \ref{lem:3.3 Christ},
\begin{align*}
&\frac{1}{\mu_p(z,\frac 1\tau)^{|\alpha|-1}}
e^{-c\frac{|z-w|^2}s} e^{-c \frac{s}{\mu_p(w,1/\tau)^2}}
\les s^{\frac 12} \frac{1}{\mu_p(z,\frac 1\tau)^{|\alpha|}}
\frac{\mu_p(w,\frac 1\tau)^{|\alpha|}}{s^{\frac 12|\alpha|}}
e^{-c\frac{|z-w|^2}s} e^{-c \frac{s}{\mu_p(w,1/\tau)^2}} \\
&\les \frac{1}{s^{\frac 12|\alpha|-\frac 12}}  \left(\frac{|z-w|}{\mu_p(w,\frac 1\tau)}\right)^M
\left(\frac{s}{|z-w|^2}\right)^{\frac M2} \left(\frac{\mu_p(w,\frac 1\tau)^2}{s}\right)^{\frac M2}
e^{-c\frac{|z-w|^2}s} e^{-c \frac{s}{\mu_p(w,1/\tau)^2}} \\
&\les \frac{1}{s^{\frac 12|\alpha| - \frac 12}}  
e^{-c\frac{|z-w|^2}s} e^{-c \frac{s}{\mu_p(w,1/\tau)^2}}.
\end{align*}
If $s<\mu_p(z,\frac 1\tau)$, the estimate is simpler because $\mu_p(z,\frac 1\tau)<|w-z|$, and we can use the Gaussian
decay term to achieve the desired bound.

To estimate $IV$, note that if
$|v-w|<|v-z|$, then $|v-z| \geq \frac 12|w-z|$, so
\begin{align}
|IV|
&\leq \int_{\atopp{|v-w|>\mu_p(w,1/\tau)}{\atopp{|v-z|>\mu_p(z,1/\tau)}{|v-w|<|v-z|}}}
|U^\gamma_w\Htp(s,v,w)  X^\alpha_z\Rtp(z,v)|\, dA(v) \nn\\
&\les \frac{1}{s^{1+\frac 12|\gamma|}}e^{-c\frac{s}{\mu_p(w,1/\tau)^2}}
\frac{e^{-c \frac{|z-w|}{\mu_p(z,1/\tau)}}} {\mu_p(z,\frac 1\tau)^{1+|\alpha|}} 
\int_{|v-w|>\mu_p(w,\frac 1\tau)} e^{-c \frac{|v-w|^2}{s}}\, dA(v) \nn\\
&\les \frac{1}{s^{1+\frac 12|\gamma|}}e^{-c\frac{s}{\mu_p(w,1/\tau)^2}}
\frac{e^{-c \frac{|z-w|}{\mu_p(z,1/\tau)}}} {\mu_p(z,\frac 1\tau)^{1+|\alpha|}} 
\int_{\mu_p(w,1/\tau)}^\infty s e^{-c \frac {r^2}s}\frac rs\, dr \nn\\
&\sim \frac{1}{s^{\frac 12|\gamma|}}e^{-c\frac{s}{\mu_p(w,1/\tau)^2}}
\frac{e^{-c \frac{|z-w|}{\mu_p(z,1/\tau)}}} {\mu_p(z,\frac 1\tau)^{1+|\alpha|}}
e^{-c \frac{\mu_p(w,1/\tau)^2}s} \label{eqn:IV est} 
\end{align}
(\ref{eqn:IV est}) easily reduces to the desired estimate.

We now have to worry about the $|z-w|\leq 2\mu_p(z,\frac 1\tau)$ case. Note that if 
$|z-w|\leq 2\mu_p(z,\frac 1\tau)$, then $\mu_p(w,\frac 1\tau) \sim \mu_p(z,\frac 1\tau)$. Similarly to above, let
$\psi_z \in \cic{\C}$ with $\supp \psi \subset D(z,4\mu_p(z,1/\tau))$, $\psi\equiv 1$ on $D(z,2\mu_p(z,1/\tau))$,
and $|D^k \psi| \leq c_k \mu_p(z,1/\tau)^{-k}$. Then
\begin{align*}
&\int_\C U^\gamma_w\Htp(s,v,w) X^\alpha_z\Rtp(z,v)\, dA(v) \\
&= \int_\C U^\gamma_w\Htp(s,v,w) X^\alpha_z\Rtp(z,v)\psi_z(v)\, dA(v) 
+\int_\C U^\gamma_w\Htp(s,v,w) X^\alpha_z\Rtp(z,v)(1-\psi_z(v))\, dA(v) \\
\end{align*}
The first term in the above equality can be estimated 
using the cancellation condition for $\Rtp$, Lemma \ref{lem:R cancel}. This
estimate is the same as the estimate of $I$.  To bound the second term,  we can split the integral up
with regions $|v-w|>|v-z|$ and $|v-w|<|v-z|$. In this manner, the estimate is very similar to 
the estimates for $III$ and $IV$.
\end{proof}

%
%
\section{Off-diagonal estimates of $\Hwtp(s,z,w)$}\label{sec:gaussian decay}

We  prove an off-diagonal
Gaussian decay estimate for $\Hwtp(s,z,w)$. On-diagonal 
Gaussian estimates, i.e.,
estimates where $|z-w| \leq \sqrt s$, are false since $\lim_{s\to\infty}\Hwtp(s,z,w) = \Ss(z,w)$. The
$1/s$ term in Gaussian decay would force $\lim_{s\to\infty}\Hwtp(s,z,w)=0$ which clearly does not occur. 

We need some notation.
Let $E^{\Boxwtp}$ be the resolution of the identity for $\Boxwtp$
(see \cite{Rudin91} or \cite{ReSi80} for details). By the spectral theorem, the operator 
\[
F(\Boxwtp) = \int_0^\infty F(\lambda)\, dE^{\Boxwtp}(\lambda).
\]
is defined so that given $f_1,f_2\in L^2(\C)$, 
$dE^{\Boxwtp}_{f_1,f_2}$ is a measure with support on
$\sigma(\Boxwtp)\subset [0,\infty)$ and 
\[
\Big( F(\Boxwtp) f_1,f_2\Big) = \int_0^\infty F(\lambda) dE^{\Boxwtp}_{f_1,f_2}.
\]
Our integral operators $F(\Boxwtp)$ can be given by integration against a kernel, and we denote the kernel
of $F(\Boxwtp)$ by $K_{F}(z,w)$. For example, $K_{e^{-s\Boxwtp}}(z,w) = \Hwtp(s,z,w)$.

The idea of the argument
is to decompose $e^{-s\Boxwtp} = F_{\ell}^0(\sqrt{s\Boxwtp}) + F_{\ell}^\infty(\sqrt{s\Boxwtp})$.
We will fix $(s,z,w)$ and choose $\ell$ and $F^0_\ell$ appropriately so that
$K_{F_\ell^0}(z,w)=0$. The decomposition is based on a finite propogation speed 
result for operators which can be
written in terms of the wave kernel $\cos(\sqrt{s\Boxwtp})$. The ``cut-off" introduces a Gaussian
decay term into $F_\ell^\infty$, and a good poinwise estimate of the tail of the Fourier transform
of a Gaussian allows us to recover the off-diagonal Gaussian decay estimate for $\Hwtp(s,z,w)$.
In addition to the tail estimate, we also need $L^2$-estimates on the integral kernel of the
resolvant $R(\ell,\Boxwtp) = (\ell I-\Boxwtp)^{-1}$, and $L^2$-estimates on the kernel of a 
composition of operators.

\subsection{Finite Propogation Speed}

The finite propogation speed of solutions to the wave equation
will follow from a local conservation of energy law.
\begin{prop}\label{prop:energy2}
Fix $(s_0,z_0) \in (0,\infty)\times\C$. Let $u$ solve the wave equation
\[
\frac{\p^2 u}{\p s^2} + \Boxwtp u =0
\]
in $\{(s,z) : 0\leq s\leq s_0 \text{ and }z\in \overline{D(z_0,s_0-s)}\}$.
If
\[
e[u](s) = \int_{D(z_0, s_0-s)} |u_s(z)|^2 + |\Zbstp u(z)|^2 \, dA(z),
\]
then $e[u](s) \leq e[u](0)$.
\end{prop}

\begin{proof} It's enough to show $\dot{e}[u](s) \leq 0$ for $0 < s < s_0$. If $dS$ is ``surface area'' measure
on the circle with outward unit normal $\nu$, we have:
\begin{align*} 
\dot{e}[u](s) &= \int_{D(z_0, s_0-s)} \frac{\p}{\p s}u_s \bar u_s
+ \frac{\p}{\p s}(\Zbstp u\overline{\Zbstp u})\, dA(z)
- \int_{\p D(z_0,s_0-s)} |u_s|^2 + |\Zbstp u|^2 \, dS \\
&= 2 \int_{D(z_0, s_0-s)} \Rre\big(u_{ss}\bar u_s  + \Zbstp u_s \overline{\Zbstp u} \big)
- \int_{\p D(z_0,s_0-s)} |u_s|^2 + |\Zbstp u|^2 \, dS \\
&= 2\Rre \left(\int_{D(z_0, s_0-s)} \hspace{-15.5pt} u_{ss}\bar u_s   + 
\Boxwtp u \bar u_s\, dA(z) +\int_{\p D(z_0,s_0-s)} \hspace{-15.5pt}u_s \overline{\Zbstp u} \nu\, dS \right)
- \int_{\p D(z_0,s_0-s)} \hspace{-18pt} |u_s|^2 + |\Zbstp u|^2 \, dS \\
&= 2\Rre\left(\int_{D(z_0, s_0-s)} \hspace{-15.5pt}\bar u_s(\underbrace{u_{ss} + \Boxwtp u}_{=0}) 
+\int_{\p D(z_0,s_0-s)}\hspace{-15.5pt} u_s \overline{\Zbstp u}\, \nu dS\right)
- \int_{\p D(z_0,s_0-s)} \hspace{-15.5pt}|u_s|^2 + |\Zbstp u|^2 \, dS \\
&= - \int_{\p D(z_0,s_0-s)} |u_s|^2 + |\Zbstp u|^2 \, dS
+ 2\int_{\p D(z_0,s_0-s)}\Rre\left( u_s \overline{\Zbstp u}\right)\, \nu dS\\
&\leq - \int_{\p D(z_0,s_0-s)} |u_s|^2 + |\Zbstp u|^2 \, dS
+ 2\int_{\p D(z_0,s_0-s)}|u_s \overline{\Zbstp u}|\, dS\\
&\leq - \int_{\p D(z_0,s_0-s)} |u_s|^2 + |\Zbstp u|^2 \, dS
+ 2\int_{\p D(z_0,s_0-s)}\frac 12 |u_s|^2 + \frac12 |\Zbstp u|^2 \, dS  \leq 0.
\end{align*}
\end{proof}
We have the following corollary.
\begin{cor}\label{cor:fps}
Fix $(s_0,z_0) \in (0,\infty)\times\C$. Let $u$ solve the wave equation
\[
\frac{\p^2 u}{\p s^2} + \Boxwtp u =0
\]
in $\{(s,z) : 0\leq s\leq s_0 \text{ and }z\in \overline{D(z_0,s_0-s)}\}$.
If $u(0,z) = \frac{\p u}{\p s}(0,z) =0$ on $D(z_0,s_0)$, then 
\[
u(s,z)\equiv0
\] 
on
the cone 
\[
\{(s,z) : 0\leq s\leq s_0: |z-z_0| < s_0-s\}.
\]
\end{cor}
\begin{proof} From Lemma \ref{prop:energy2},
\[
\int_{D(z_0,s_0-s)} |u_s(s,z)|^2 + |\Zbstp u(s,z)|^2\, dA(z)
\leq \int_{D(z_0,s_0)}|u_s(0,z)|^2 + |\Zbstp u(0,z)|^2\, dA(z).
\]
However, $u_s(0,z)=0$ on $D(z_0,s_0)$. Also, since $u(0,z)=0$ in $D(z_0,s_0)$, $\Zbstp u(0,z)=0$
in $D(z_0,s_0-s)$. The corollary follows quickly.
\end{proof}

Let $f$ be a suitably nice function defined on $\C$. Use the spectral theorem to define
$\Vtp f(s,z) = \cos(s\sqrt{\Boxwtp})[f](z)$. $\Vtp$ is a contraction on 
$L^2(\C)$, so we can write
\[
\Vtp[f](s,z) = \int_{\C} \Vtp(s,z,w) f(w)\, dA(w)
\]
for some distributional kernel $\Vtp(s,z,w)$. In our new notation, Corollary
\ref{cor:fps} becomes
\begin{cor}\label{cor:fps2}
If $f(z)=0$ on $D(z_0,s_0)$, then
\[
\Vtp f(s,z) \equiv 0
\]
on the cone
\[
\{(s,z) : 0\leq s\leq s_0: |z-z_0| < s_0-s\}.
\]
\end{cor}
The finite propagation speed of solutions to the wave equation implies a certain support
of the wave kernel $\Vtp(s,z,w)$.
\begin{cor}\label{cor:W support} 
\[
 \Vtp(s,z,w)=0 \text{ if } |z-w|>s.
\]
\end{cor}

\begin{proof} Fix $s\in [0,\infty)$ and $z\in\C$.
Let $f\in C^\infty_c(\C)$ be defined so that $\supp f \subset \C\setminus D(z,s)$. Then with $s_0=s$ and
$z_0=z$, 
Corollary \ref{cor:fps2} shows $\Vtp f(s,z)=0$. As a result
\begin{align*}
0&=\Vtp f(s,z) = \int_{D(z,s)} \Vtp(s,z,w) f(w) \,dA(w) + \int_{\C\setminus D(z,s)} \Vtp(s,z,w) f(w)\, dA(w)\\
&= \int_{\C\setminus D(z,s)} \Vtp(s,z,w) f(w)\, dA(w).
\end{align*}
Since $f$ was arbitrary, $\Vtp(s,z,w)=0$ on $\{w\in \C\setminus D(z,s)\}$.
\end{proof}

The finite propogation speed of the wave kernel allows us to adapt Sikora's technique
\cite{Sik04} for off-diagonal Gaussian estimates.
\subsection{Spectral Theorem Facts}
Now we collect some useful results. The support condition of the wave kernel has the following consequence which
is a crucial tool for our decomposition of $e^{-s\Boxwtp}$.
\begin{prop}\label{prop:support condition}
Let $F$ be a bounded, even function defined on $\R$  whose Fourier transform satisfies 
$\supp \hat F \subset [-L,L]$. Then $K_{F(\sqrt{s\Boxwtp})}(z,w)=0$ if $|z-w|>s^{1/2}L$.
\end{prop}

\begin{proof}Since $F$ is even, the Fourier inversion formula can be written
$F(\sqrt x) = \frac{1}{2\pi}\int_{\R}\hat F(t)\cos(t\sqrt x)\, dt$, and the spectral theorem implies
\[
F(\sqrt {s\Boxwtp}) = \frac{1}{2\pi}\int_{\R}\hat F(t)\cos(t\sqrt{s\Boxwtp})\, dt.
\]
Thus, for $f_1,f_2\in C^\infty_c(\C)$,
\begin{equation}\label{eqn:support}
\Big( F(\sqrt{s\Boxwtp}) f_1,f_2\Big) 
=\int_\C\int_\C	 \underbrace{\int_{-L}^L \hat F(t) \Vtp(s^{1/2}t,z,w)\, dt}_{K_{F(\sqrt{s\Boxwtp})}} f_1(w)
\overline{f_2(z)}\, dA(w) dA(z).
\end{equation}
Proposition \ref{prop:support condition} follows quickly from  \eqref{eqn:support} and
Corollary \ref{cor:W support}.
\end{proof} 

Next we investigate products.
\begin{prop}\label{prop:products} If $F$ and $G$ are bounded functions on $\R$,
\begin{equation}\label{eqn:products}
\| K_{FG}(z,\cdot) \|_{L^2(\C)} \leq \|G\|_{L^\infty(\C)} \|K_{F}(z,\cdot)\|_{L^2(\C)}.
\end{equation}
\end{prop}

\begin{proof} Since $G$ is bounded, $G(\Boxwtp) = \int_{0}^\infty G(\lambda)\, dE^{\Boxwtp}(\lambda)$, and
it is clear that $\|G(\Boxwtp)\|_{L^2\to L^2} \leq \|G\|_{L^\infty(\C)}$. Also, fixing $z,w\in\C$ and
writing $K_{F}(z,v) = K^z_F(v)$ to emphasize that $z$ is fixed, we have
\[
K_{FG}(z,w) = \int_{\C}K_F(z,v)K_G(v,w)\, dA(v) = \int_{\C} K^z_F(v) K_G(v,w)\, dA(v) 
= G^*(\sqrt{\Boxwtp})[K_F^z](w).
\]
But $\|G\|_{L^2\to L^2} = \|G^*\|_{L^2\to L^2}$, so
\[
\| K_{FG}(z,\cdot)\|_{L^2(\C)} = \|G^*(\sqrt{\Box})[K_F^z] \|_{L^2(\C)} \leq 
\|G(\sqrt{\Boxwtp})\|_{L^2\to L^2}\|K_F(z,\cdot)\|_{L^2(\C)}.
\]
\end{proof}

\subsection{$L^2$-estimates of the resolvant and a decay estimate}

As mentioned above, we need a decay result for the Fourier transform on the tail of a Gaussian in one dimension.
Let $\ell\geq 1$ and $\vp_\ell\in C^\infty(\R)$ so that $\vp_\ell$ is even, $\vp_\ell(x) =1$ if $|x|\geq\ell$,
and $\vp(x)=0$ if $|x|\leq \frac 12\ell$. Also, let $0\leq \vp_\ell \leq 1$ and $|\vp^{(n)}| \leq a_n \ell^{-n}$.
\begin{prop}\label{prop:cut off gaussian} For each $N\in\N$, there exists a constant $C_N$ so that if $\lam\geq 0$:
\[
\left| \mathcal{F}^{-1}[\vp_\ell e^{-\frac{x^2}4}](\lam)\right| = 
\left| \frac 1{2\pi}\int_{\R}\vp_\ell(\xi) e^{-\frac{\xi^2}4} e^{i\lam \xi}\, d\xi \right|
\leq C_N e^{-\frac{\ell^2}{16}} \frac{1}{(\ell^2+\lam^2)^{N/2}}.
\]
\end{prop}

\begin{proof} The proof follows from an integration by parts. It is enough to estimate the integral
on $[0,\infty)$. Integrating by parts $N$ times yields:
\begin{align*}
\int_0^\infty \vp_\ell(\xi) e^{-\frac{\xi^2}4+i\lam \xi}\, d\xi
&= \int_0^\infty \left( \frac{1}{\xi/2 - i\lam}
\left(\cdots \left( \frac{1}{\xi/2-i\lam}\vp_\ell(\xi)\right)'\cdots\right)'\right)'
e^{-\frac{\xi^2}4 + i\lam \xi}\, d\xi \\
&= \sum_{j=0}^N c_j \int_0^\infty 
\frac{1}{(\xi/2-i\lam)^{N+j}}\vp_\ell^{(N-j)}(\xi)e^{-\frac{\xi^2}4 + i\lam \xi}\, d\xi. \\
\end{align*}
If $j=N$,
\begin{align*}
\Big| \int_{\frac\ell2}^\infty (\xi/2 - i\lam)^{-2N}\vp_\ell(\xi)  e^{-\frac{\xi^2}4 + i\lam \xi}\, d\xi\Big|
&\leq e^{-\frac{\ell^2}{16}} \int_{\ell}^\infty \Big|\frac \xi2 - i\lam\Big|^{-2N}\, d\xi \\
&\leq c_N e^{-\frac{\ell^2}{16}} \int_{\ell}^\infty (\xi + \lam)^{-2N}\, d\xi \\
&\leq c_N e^{-\frac{\ell^2}{16}}(\ell+\lam)^{-2N+1} = C_N e^{-\frac{\ell^2}{16}}(\ell^2+\lam^2)^{-N+\frac12}.
\end{align*}
If $j<N$,
\begin{align*}
\Big|\int_0^\infty \frac{1}{(\xi/2-i\lam)^{N+j}}\vp_\ell^{(N-j)}(\xi)e^{-\frac{\xi^2}{4} 
+ i\lam \xi}\, d\xi\Big|
&\leq c_N e^{-\frac{\ell^2}{16}}\ell |\ell - i\lam|^{-j-N}\ell^{-N+j} \\
&\leq c_N e^{-\frac{\ell^2}{16}}(\ell^2+ \lam^2)^{-N/2-j/2}\ell^{-N+j+1}.
\end{align*}
The $j=0$ case is the worst estimate from the previous two inequalities, and we see that
\[
\left| \mathcal{F}^{-1}[\vp_\ell e^{-\frac{x^2}{16}}](\lam)\right| 
\leq C_N e^{-\frac{\ell^2}{16}} (\ell^2+\lam^2)^{-N/2}\ell^{-N+1}.
\]
\end{proof}

Next, we need $L^2$-bounds on the kernel of the resolvent 
$R(\lam,\Boxwtp) = (\lam I + \Boxwtp)^{-1}$.  $\Boxwtp$ has nonnegative eigenvalues, so 
$R(\lam,\Boxwtp)$ is a bounded, well-defined operator on $L^2(\C)$ if $\lam >0$. By the spectral theorem,
if $\lam>0$,
\begin{equation}\label{eqn:resolvent}
R(\lam,\Boxwtp)[f] = \lim_{\ep\to0} \int_{\ep}^{1/\ep} (e^{-\lam s}e^{-s\Boxwtp})[f]\, ds
= \lim_{\ep\to0} \int_\ep^{1/\ep} e^{-\lam s} H_s[f]\, ds.
\end{equation}
If $r_\lam(z,w)$ is the integral kernel of $R(\lam,\Boxwtp)$, then it follows from \eqref{eqn:resolvent} that
\begin{equation}\label{eqn:heat resolvent}
r_\lam(z,w) = \int_0^\infty e^{-\lam s} \Hwtp(s,z,w)\, ds.
\end{equation}
We will need to compute $\int_\C |r_\lam(z,w)|^2\, dA(w)$. A first step in that computation is
\begin{prop}\label{prop:L2 heat} 
Let $Y^{\alpha}$ be a product of operators $Y = X_j$, $j=1,2$ if acting in $z$ or
$U_j$, $j=1,2$ if acting in $w$. Then there exists a constant $C_{|\alpha|}$ so that
\[
\| Y^\alpha \Hwtp(s,z,\cdot)\|_{L^2(\C)} \leq C_{|\alpha|} \max\{ s^{-\frac 12-\frac 12|\alpha|},
 \mu_p(z,1/\tau)^{-1-|\alpha|}\}.
\]
\end{prop}

\begin{proof} Recall that $\Hwtp(s,z,w) = \Gwtp(s,z,w) + \Ss(z,w)$. By Corollary \ref{cor:R,S decay}
and the change of variables that
$x = c(z-w)/\mu(z,1/\tau)$, we have
\begin{align}
\int_{\C}|Y^\alpha\Ss(z,w)|^2 \, dA(w) 
&\les \frac{1}{\mu(z,1/\tau)^{4+2|\alpha|}} \int_{\C} e^{-c\frac{|z-w|}{\mu(z,1/\tau)}}\, dA(w)\nn\\
&= \frac{1}{\mu(z,1/\tau)^{2+2|\alpha|}}\int_{\C} e^{-|x|}\,dA(x) \sim \frac{1}{\mu(z,1/\tau)^{2+2|\alpha|}}.
\label{eqn:L^2 Szego}
\end{align}
Also, from Theorem \ref{thm:G wig and deriv}, we see that one of the bounds on $\Gwtp(s,z,w)$ is (up to constants) the
same bound for $\Ss(z,w)$ and hence already estimated. 
Also, integrating
the Gaussian decay bound of $Y^\alpha \Gwtp(s,z,w)$, we have (with $x =\sqrt{c} (z-w)/\sqrt s$)
\begin{equation}\label{eqn:Gw s bound}
\int_{\C} \frac{1}{s^{2+|\alpha|}}e^{-c|z-w|^2/s}\, dA(w) 
= \frac 1{s^{1+|\alpha|}} \int_{\C}e^{-|x|^2}\, dA(x) \sim 
\frac 1{s^{1+|\alpha|}}.
\end{equation}
Taking the square roots of the estimates \eqref{eqn:L^2 Szego} 
and \eqref{eqn:Gw s bound} finishes the proof.
\end{proof}

Using Proposition \ref{prop:L2 heat}, we can find an $L^2$ bound on $r_\lam(z,w)$.
\begin{prop}\label{prop:L2 rlam} There exists a constant $C$ so that
\[
\|r_{\lam}(z,\cdot)\|_{L^2(\C)} \leq C\max\left\{\frac{1}{\lam\mu(z,1/\tau)},\frac 1{\lam^{1/2}}\right\}.
\]
\end{prop}

\begin{proof} From \eqref{eqn:heat resolvent} and Proposition \ref{prop:L2 heat}, it follows from Minkowski's Inequality for Integrals that
\begin{align*}
\left(\int_{\C} |r_\lam(z,w)|^2\, dA(w)\right)^{\frac 12}
& = \left(\int_{\C}\left|\int_0^\infty e^{-\lam s} \Hwtp(s,z,w)\, ds\right|^2\, dA(w)\right)^{\frac 12} \\
&\leq \int_0^\infty e^{-\lam s} \|\Hwtp(s,z,\cdot)\|\, ds 
\leq \int_0^\infty e^{-\lam s}\max\{s^{-1/2},\mu_p(z,1/\tau)^{-1}\}\,ds.
\end{align*}
Certainly $\int_0^\infty e^{-\lam s}\mu(z,1/\tau)^{-1}\, ds = \lam^{-1}\mu(z,1/\tau)^{-1}$. 
Also, 
\[
\int_0^\infty s^{-1/2}e^{-\lam s} \leq \int_0^{1/\lam} s^{-1/2}\, ds + \int_{1/\lam}^\infty \lam^{1/2} 
e^{-\lam s}\, ds \sim \lam^{-1/2}.
\] 
\end{proof}

\subsection{Relationship between the heat and wave kernels and a 
decomposition of the heat kernel}
The next part of our discussion explores the relationship between the heat kernel
$\Hwtp(s,z,w)$ and the wave kernel $\Vtp(s,z,w)$. If $f\in L^1(\R)$,
\[
\hat f(\xi) = \int_{\R} e^{-i \xi x}f(x)\, dx.
\]
From this formula, it quickly follows that
\begin{equation}\label{eqn:four trans}
\frac{1}{2\sqrt{\pi s}}\int_0^\infty e^{-r^2/4s} \cos r\xi\, dr = e^{-s\xi^2}.
\end{equation}
\eqref{eqn:four trans} is the basis for the relationship between the wave kernel and the heat kernel. Since
$\Boxwtp$ is self-adjoint and positive (the latter condition meaning $\big(\Boxwtp f,f\big)\geq 0$ for all 
$f\in\Dom(\Boxwtp)$), 
the spectral theorem implies that $\sqrt{\Boxwtp}$ is a well-defined operator. Even further,
the spectral theorem and \eqref{eqn:four trans} prove that
\begin{equation}\label{eqn:he wa}
\frac{1}{2\sqrt{\pi s}} \int_0^\infty \Big(\cos(r\sqrt{\Boxwtp}) f_1,f_2\Big) e^{-r^2/4s}\, dr
= \Big( e^{-s\Boxwtp} f_1,f_2\Big).
\end{equation}
From the kernel side, equation \eqref{eqn:he wa} becomes
\begin{multline}\label{eqn:he wa ker}
\frac{1}{2\sqrt{\pi s}} \int_0^\infty\int_{\C}\int_{\C} \Vtp(r,z,w) f_1(w)\overline{f_2(z)}
e^{-r^2/4s}\,dA(w)dA(z)dr \\
 = \int_{\C}\int_{\C} \Hwtp(s,z,w) f_1(w)\overline{f_2(z)}\, dA(w) dA(z).
\end{multline} 
Thus, from \eqref{eqn:he wa ker} and the support condition from Corollary \ref{cor:W support}, it follows that
\begin{prop}\label{prop:he wa ker}
\begin{equation}\label{eqn:he wa ker2}
\Hwtp(s,z,w) = \frac{1}{2\sqrt{\pi s}} \int_{0}^\infty \Vtp(r,z,w) e^{-r^2/4s}\, dr
= \frac{1}{2\sqrt{\pi s}} \int_{|z-w|}^\infty \Vtp(r,z,w) e^{-r^2/4s}\, dr.
\end{equation}
\end{prop}

We are finally ready to decompose $e^{-s\Boxwtp}$. Let $\ell>1$ (to be chosen later) and let
\[
\widehat{F^\infty_\ell}(\xi) = \sqrt\pi \vp_{\ell}(\xi)e^{-\frac{\xi^2}4}
\qquad\text{and}\qquad 
\widehat{F^0_\ell}(\xi) = \sqrt\pi (1-\vp_{\ell}(\xi))e^{-\frac{\xi^2}4}.
\]
Then $F^0_\ell$ and $F^\infty_\ell$ are bounded, even functions and
\[
e^{-s\Boxwtp} = F^0_\ell(\sqrt{s\Boxwtp}) + F^\infty_\ell(\sqrt{s\Boxwtp}).
\]
Next, $\supp(1-\vp_\ell) \subset (-\ell,\ell)$, so by Proposition \ref{prop:support condition},
\begin{equation}\label{eqn:good support}
\supp K_{F^0_\ell(\sqrt{s\Boxwtp})}(\cdot,\cdot) \subset \{(z,w) : |z-w| < s^{1/2}\ell\}.
\end{equation}

Fix $z,w\in\C$ and $s>0$ so that $|z-w|>s^{\frac 12}$.
Set $\ell = \ell(z,w) = |z-w|/s^{1/2}$. By \eqref{eqn:good support},
\[
K_{F^0_\ell(\sqrt{s\Boxwtp})}(z,w)=0,
\]
so
\begin{equation}\label{eqn:reduction}
\Hwtp(s,z,w) = K_{F^\infty_\ell(\sqrt{s\Boxwtp})}(z,w).
\end{equation}

For simplicity, set $F(x) = F^\infty_\ell(x)$. Note that Proposition 
\ref{prop:cut off gaussian} applies to $F(x)$.
For the next step in the argument, we would like to
take $\sqrt{F}$, but that is unlikely to exist. While $F$ is $\R$-valued, there is no reason to believe that
$F\geq 0$. Instead, we can write $F(\xi) = |F(\xi)|e^{i\theta(\xi)}$, and our replacement for the ``square root'' of 
$F$ is $f(\xi) = |F(\xi)|^{1/2} e^{i\theta(\xi)/2}$ with our branch cut taken so that $0\leq\theta<2\pi$.
In fact, $\theta(\xi) = 0$ or $\pi$. 
From Proposition \ref{prop:cut off gaussian}, 
we have (taking a square root of the estimate) that for all $N\geq 1$
\begin{equation}\label{eqn:square root}
\sup_{\lam\geq 0} |f(\lam)(\ell^2+\lam^2)^{N/4}| \leq C_N e^{-c \ell^2} =C_N e^{-c \frac{|z-w|^2}s}.
\end{equation}

Next, we would like to use Proposition \ref{prop:products} 
and (\ref{eqn:square root}), our weighted $L^\infty$-estimate of $f$,  to bound
$\Hwtp(s,z,w)$. We have:
\begin{equation}\label{eqn:HK relation}
|\Hwtp(s,z,w)| = |K_{F(\sqrt{s\Boxwtp})}(z,w)| 
\leq \|K_{f(\sqrt{s\Boxwtp})}(z,\cdot)\|_{L^2(\C)} \|K_{f(\sqrt{s\Boxwtp})}(\cdot,w)\|_{L^2(\C)}
\end{equation} 
Evaluating each term separately, we have (with $\lam = \sqrt{s\Boxwtp}$ and $N=4$) by
Proposition \ref{prop:products}
\begin{align}
\| K_{f(\sqrt{s\Boxwtp})} (z,\cdot)\|_{L^2(\C)}
&= \|K_{f(\sqrt{s\Boxwtp})(\ell^2 I + s\Boxwtp)(\ell^2 I+s\Boxwtp)^{-1}}(z,\cdot)\|_{L^2(\C)}\nn\\
&\les \sup_{\lam>0}\left|f(\lam)\left(\ell^2+\lam^2\right)\right|
\|K_{(\ell^2 + s\Boxwtp)^{-1}}(z,\cdot)\|_{L^2(\C)}. \label{eqn:L2 Kf}
\end{align}
Next, since $\ell = \frac{|z-w|}{s^{1/2}}>1$, from Proposition \ref{prop:L2 rlam} we have 
\begin{align}
\|K_{(\ell^2I + s\Boxwtp)^{-1}}(z,\cdot)\|_{L^2(\C)}
&= \|K_{ \frac1s(\frac{\ell^2}{s}I + \Boxwtp)^{-1}}(z,\cdot)\|_{L^2(\C)}
= \frac{1}{s} \|r_{\frac{\ell^2}{s}}(z,\cdot)\|_{L^2(\C)}\nn\\
&\les \frac{1}{s} \max\left\{\frac{s}{\ell^2 \mu(z,1/\tau)},\frac{s^{1/2}}{\ell}\right\}
\les \max\left\{\frac{1}{\mu(z,1/\tau)}, \frac1{s^{1/2}}\right\}. \label{eqn:L2 K}
\end{align}
Putting our estimates together finally allows us to prove 
%
%
\begin{thm}\label{thm:off diag est} There exists constants $c, C>0$ so that if $|z-w|>s^{1/2}$,
\[
\big| \Hwtp(s,z,w) \big| \leq C e^{-c \frac{|z-w|^2}{s}}
\max\left\{ \frac 1{\mu_p(z,1/\tau)\mu_p(w,1/\tau)}, \frac{1}{s}\right\}.
\]
\end{thm}

For a discussion on a quantitative estimate of $c$, see \cite{Sik04}.
\begin{proof} The result follows immediately from combining
\eqref{eqn:HK relation}, \eqref{eqn:square root}, \eqref{eqn:L2 Kf}, and \eqref{eqn:L2 K}.
\end{proof}
\section{Proof of Theorem \ref{thm:H wig and deriv}}
We first prove off-diagonal estimates for $X^\alpha_z U^\gamma_w \Hwtp(s,z,w)$.
We need the following result.
\begin{lem}\label{lem:mult maxes}
Let $a, b >0$, $|z-w| \geq \mu_p(z,1/\tau)$, and $|z-w| \geq\mu_p(w,1/\tau)$.
Then there exists a constant $C_{a,b}$ so that
\begin{multline*}
e^{-c \frac{|z-w|^2}{s}} \max \left\{ \frac{1}{s^{a}}, \frac{1}{\mu_p(z, 1/\tau)^{2a}} \right\}
\cdot \max\left\{\frac{e^{-c\frac{s}{\mu_p(z,1/\tau)^2}}}{s^{b}}, 
\frac{e^{-c\frac{|z-w|}{\mu_p(z,1/\tau)}}}{\mu_p(z, 1/\tau)^{2b}} \right\} \\
\le C_{a,b} e^{-c \frac{|z-w|^2}{s}}\max\left\{\frac{e^{-c\frac{s}{\mu_p(z,1/\tau)^2}}}{s^{a+b}}, 
\frac{e^{-c\frac{|z-w|}{\mu_p(z,1/\tau)}}}{\mu_p(z, 1/\tau)^{2(a+b)}} \right\}.
\end{multline*}
\end{lem}

\begin{proof}
If $s^{-1}\geq\mu_p(z,1/\tau)^{-2}$ and $
\frac{e^{-c\frac{s}{\mu_p(z,1/\tau)^2}}}{s^b}\geq 
\frac{e^{-c\frac{|z-w|}{\mu_p(z,1/\tau)}}}{\mu_p(z, 1/\tau)^{2b}}$
or $s^{-1}\leq\mu_p(z,1/\tau)^{-2}$ and 
$\frac{e^{-c\frac{s}{\mu_p(z,1/\tau)^2}}}{s^b}\leq
\frac{e^{-c\frac{|z-w|}{\mu_p(z,1/\tau)}}}{\mu_p(z, 1/\tau)^{2b}}$,
then there is nothing to prove. Thus we have two cases to consider. Assume that
$s^{-1}\geq\mu_p(z,1/\tau)^{-2}$ and 
$\frac{e^{-c\frac{s}{\mu_p(z,1/\tau)^2}}}{s^b}\leq
\frac{e^{-c\frac{|z-w|}{\mu_p(z,1/\tau)}}}{\mu_p(z, 1/\tau)^{2b}}$. Note that with
a decrease in $c$, $e^{-c \frac {|z-w|^2}{s}} \leq
\frac{s^{a}}{|z-w|^{2a}}|z-w| e^{-c \frac {|z-w|^2}{s}}$. Since $|z-w|\geq \mu_p(z,1/\tau)$,
\[
e^{-c \frac {|z-w|^2}s} \frac 1{s^a} 
\frac{e^{-c\frac{|z-w|}{\mu_p(z,1/\tau)}}}{\mu_p(z, 1/\tau)^{2b}}
\leq 
\frac{s^{a}}{|z-w|^{2a}} e^{-c \frac {|z-w|^2}{2s}}
\frac 1{s^{a}} \frac{e^{-c\frac{|z-w|}{\mu_p(z,1/\tau)}}}{\mu_p(z, 1/\tau)^{2b}}.
\]
Now assume that $s^{-1}\leq\mu_p(z,1/\tau)^{-2}$ and 
$\frac{e^{-c\frac{s}{\mu_p(z,1/\tau)^2}}}{s^b}\geq
\frac{e^{-c\frac{|z-w|}{\mu_p(z,1/\tau)}}}{\mu_p(z, 1/\tau)^{2b}}$. 
The argument is similar to the first case, except we use the fact that 
$ e^{-c\frac{s}{\mu_p(z,1/\tau)^2}} \les \frac{\mu_p(z,1/\tau)^{2a}}{s^{a}}
e^{-c\frac{s}{\mu_p(z,1/\tau)^2}}$ with a decrease in the constant $c$.
\end{proof}

The bound on $X^\alpha_z U^\gamma_w \Hwtp(s,z,w)$ is proven in two stages. We first prove:
\begin{lem}\label{lem:Hwtp deriv est}
Let $z, w\in\C$ so that $|z-w| \geq 2\mu_p(z, 1/\tau)$ and 
$|z-w| \geq2\mu_p(w,1/\tau)$. Let
$X^\alpha_z$ be a product of $|\alpha|$ operators of the form $X=X_j$, $j=1,2$ acting in the $z$
variable. Then there exists a constant $C = C_{|\alpha|}$ so that
\[
|X^\alpha_z \Hwtp(s,z,w)| \leq C e^{-c\frac{|z-w|^2}s}
\max \left\{ \frac{e^{-c \frac{s}{\mu_p(z,1/\tau)^2}}}{ s^{1+\frac 12|\alpha|}},
\frac{e^{-c\frac{|z-w|}{\mu_p(z,1/\tau)^2}}}{\mu_p(z,1/\tau)^{2+|\alpha|}} \right\}
\]
\end{lem}
\begin{proof}
The on-diagonal result for $|\alpha|=0$  follows quickly from Theorem \ref{thm:G wig and deriv},
Corollary \ref{cor:R,S decay}, and Lemma \ref{lem:mult maxes}. The off-diagonal result follows 
immediately from Theorem \ref{thm:off diag est}, Theorem \ref{thm:G wig and deriv},
Corollary \ref{cor:R,S decay}, and Lemma \ref{lem:mult maxes}. 

We turn our attention to the case $|\alpha|\geq 1$.
Let $\delta = \min\{s^{1/2},\mu_p(z,1/\tau)\}$. 
Let $\vp_z\in  \cic\C$ so that $\vp_z \equiv 1$ on $D(z,\delta)$,
$\vp\equiv 0$ off of $D(z,\frac 32 \delta)$, and 
$|D^n \vp_z| \leq c_n \delta^{-n}$. 
By Theorem \ref{thm:Poincare},
\begin{align*}
|X^\alpha_z& \Hwtp(s,z,w)|
\les \frac{1}{\delta} \sum_{|I|\leq 2} \delta^{|I|}
\| X^I_z [X^\alpha_z\Hwtp(s,\cdot,w)\vp_z] \|_{L^2(\C)}\\
&\leq \frac{1}{\delta} \sum_{|I|\leq 2}\sum_{|J|\leq |I|} \delta^{|I|}
\| X^J_zX^\alpha_z\Hwtp(s,\cdot,w) D^{|I|-|J|}\vp_z \|_{L^2(\C)} \\
&= \frac{1}{\delta} \sum_{|I|\leq 2}\sum_{|J|\leq |I|} \delta^{|I|} 
\Big|\Big( X^\alpha_z X^J_z\big( X^J_z X^\alpha_z\Hwtp(s,\cdot,w) D^{|I|-|J|}\vp_z\overline{D^{|I|-|J|}\vp_z}\big),
\Hwtp(s,\cdot,w) \Big)\Big|^{\frac 12} \\
&\leq  \frac{1}{\delta} \hspace{-31pt}
\sum_{\atopp{|I|\leq 2}{\atopp{|J|\leq |I|}{|J_1|+|J_2|+|J_3|=|J|+|\alpha|}}} 
\hspace{-31pt} \delta^{|I|} 
\|X^{J_1}_z X^J_z X^\alpha_z\Hwtp(s,\cdot,w) D^{|I|-|J|+|J_2|}\vp_z D^{|I|-|J|+|J_3|}\vp_z
\|_{L^2(\C)}^{\frac 12}
\|\Hwtp(s,\cdot,w)\|_{L^2(\supp\vp_z)}^{\frac 12}\\
&\leq \frac{1}{\delta} \hspace{-31pt}
\sum_{\atopp{|I|\leq 2}{\atopp{|J|\leq |I|}{|J_1|+|J_2|+|J_3|=|J|+|\alpha|}}} 
\hspace{-31pt} \delta^{|I|} \frac{1}{\delta^{\frac 12(|I|-|J|+|J_2|+|I|-|J|+|J_3|)}}
\|X^{J_1}_z X^J_z X^\alpha_z\Hwtp(s,\cdot,w)\|_{L^2(\C)}^{\frac 12}
\|\Hwtp(s,\cdot,w)\|_{L^2(\supp\vp_z)}^{\frac 12}.
\end{align*}
The choice of $\delta$ is small enough so that 
\[
\|\Hwtp(s,\cdot,w) \|_{L^2(\supp\vp_z)}
\les \delta\sup_{\zeta \in D(z,\frac 32\delta)} |\Hwtp(s,\zeta,w)|
\sim \delta|\Hwtp(s,z,w)|
\]
with a decrease in the constants $c_j$.
The estimate from the other terms follows from Proposition \ref{prop:L2 heat}, and we see
\begin{align*}
\|X^{J_1}_z X^J_z X^\alpha_z\Hwtp(s,\cdot,w)\|_{L^2(\C)}
&\leq
\max\left\{\frac{1}{s^{\frac 12(1+|\alpha|+|J|+|J_1|)}},\frac{1}{\mu_p(z,1/\tau)^{1+|\alpha|+|J|+|J_1|}}\right\}\\
&= \delta^{-1-|\alpha|-|J|-|J_1|}.
\end{align*}
Combining our previous three estimates and using the fact we have already
proved the result for $\alpha=0$, we have
\begin{align*}
&|X^\alpha_z \Hwtp(s,z,w)|\\
&\leq \frac{1}{\delta} 
\hspace{-31pt}\sum_{\atopp{|I|\leq 2}{\atopp{|J|\leq |I|}{|J_1|+|J_2|+|J_3|=|J|+|\alpha|}}}\hspace{-31pt}
\delta^{|I|} 
\delta^{-|I|+|J|-\frac12(|J_2|+|J_3|+1+|\alpha|+|J|+|J_1|)}
\delta^{\frac 12}e^{-c\frac{|z-w|^2}s}
\max \left\{ \frac{e^{-c \frac{s}{\mu_p(z,1/\tau)^2}}}{ s^{\frac 12}},
\frac{e^{-c\frac{|z-w|}{\mu_p(z,1/\tau)^2}}}{\mu_p(z,1/\tau)} \right\}\\
&\les \frac{1}{\delta^{1+|\alpha|}}
e^{-c\frac{|z-w|^2}s}
\max \left\{ \frac{e^{-c \frac{s}{\mu_p(z,1/\tau)^2}}}{ s^{\frac 12}},
\frac{e^{-c\frac{|z-w|}{\mu_p(z,1/\tau)^2}}}{\mu_p(z,1/\tau)} \right\}.
\end{align*}
Since $\frac1\delta = \max \{ \frac{1}{s^{1/2}}, \frac{1}{\mu_p(z,1/\tau)}\}$, the estimate simplifies to
\begin{align*}
|X^\alpha_z \Hwtp(s,z,w)|
&\les \max\left\{ \frac{1}{s^{1/2 + \frac 12|\alpha|}}, \frac{1}{\mu_p(z,1/\tau)^{1+|\alpha|}}\right\}
e^{-c\frac{|z-w|^2}s}
\max \left\{ \frac{e^{-c \frac{s}{\mu_p(z,1/\tau)^2}}}{ s^{\frac 12}},
\frac{e^{-c\frac{|z-w|}{\mu_p(z,1/\tau)^2}}}{\mu_p(z,1/\tau)}  \right\}.
\end{align*}
From Lemma \ref{lem:mult maxes}, it follows that
\[
|X^\alpha_z \Hwtp(s,z,w)|
\les e^{-c\frac{|z-w|^2}s}
\max \left\{ \frac{e^{-c \frac{s}{\mu_p(z,1/\tau)^2}}}{ s^{1 + \frac 12|\alpha|}},
\frac{e^{-c\frac{|z-w|}{\mu_p(z,1/\tau)^2}}}{\mu_p(z,1/\tau)^{2+|\alpha|}} \right\}.
\]
\end{proof}

\begin{proof}[Proof (Theorem \ref{thm:H wig and deriv}).] Since $\Hwtp(s,z,w) = \Gwtp(s,z,w) + \Ss(z,w)$, 
the on-diagonal estimates ($|z-w|\leq 2\mu_p(z, 1/\tau)$) follow from the bounds of
the Szeg\"o kernel and Theorem \ref{thm:G wig and deriv}. 

By Theorem \ref{thm:Gw Hw special derivatives}, the only case left to consider is the bound for
$|Y^\alpha \Hwtp(s,z,w)|$.
We can write $Y^\alpha = X^\beta_w U^\gamma_w$ for some multi-indices $\beta$ and $\gamma$.
Also, from Theorem \ref{lem:Hwtp deriv est} and the fact that $\Hwtp(s,z,w) = \overline{\Hwtp(s,w,z)}$,
it follows that we have the desired estimate for $|U^\gamma_w \Hwtp(s,z,w)| = |\overline{X^\gamma_w \Hwtp(s,w,z)}|$.
Repeating the argument of Lemma \ref{lem:Hwtp deriv est} with $U^\gamma_w \Hwtp(s,z,w)$ taking the place
of $\Hwtp(s,z,w)$ finishes the proof.
\end{proof}

\bibliographystyle{alpha}
\bibliography{mybib}

\end{document}